\documentclass[12pt]{amsart}
\usepackage{amssymb,amsmath,mathrsfs}
\usepackage[ps2pdf,colorlinks=true,urlcolor=blue,
citecolor=red,linkcolor=blue,linktocpage,pdfpagelabels,bookmarksnumbered,bookmarksopen]{hyperref}
\usepackage{epsfig,graphicx,color,mathrsfs}
\usepackage[english]{babel}
\usepackage[left=2.7cm,right=2.7cm,top=3.1cm,bottom=3.1cm]{geometry}

\newcommand{\N}{{\mathbb N}}

\newcommand{\R}{{\mathbb R}}

\renewcommand{\a }{\alpha}

\renewcommand{\b }{\beta }

\renewcommand{\d }{\delta }

\newcommand{\eps}{\varepsilon}

\renewcommand{\t}{\theta}
\renewcommand{\O}{\Omega}

\numberwithin{equation}{section}
\newtheorem{theorem}{Theorem}[section]
\newtheorem{proposition}[theorem]{Proposition}
\newtheorem{lemma}[theorem]{Lemma}

\newtheorem{definition}[theorem]{Definition}
\newtheorem{remark}[theorem]{Remark}
\theoremstyle{definition}

\newcommand{\brm}{\begin{remark}\rm}
\newcommand{\erm}{\end{remark}}
\newcommand{\brms}{\begin{remark}\rm}
\newcommand{\erms}{\end{remark}}
\newcommand{\bte}{\begin{theorem}}
\newcommand{\ete}{\end{theorem}}
\newcommand{\bpr}{\begin{proposition}}
\newcommand{\epr}{\end{proposition}}
\newcommand{\ble}{\begin{lemma}}
\newcommand{\ele}{\end{lemma}}
\newcommand{\beq}{\begin{equation}}
\newcommand{\eeq}{\end{equation}}
\newcommand{\bdm}{\begin{displaymath}}
\newcommand{\edm}{\end{displaymath}}
\numberwithin{equation}{section}

\newcommand{\bos}{\begin{remark}\rm}
\newcommand{\eos}{\end{remark}}

\newcommand{\ben}{\begin{enumerate}}
\newcommand{\een}{\end{enumerate}}

\renewcommand{\a }{\alpha }
\renewcommand{\b }{\beta }
\renewcommand{\d}{\delta }

\newcommand{\e }{\varepsilon }
\newcommand{\g }{\gamma}

\newcommand{\n }{\nabla }

\newcommand{\s }{\sigma }

\renewcommand{\t }{\tau }

\renewcommand{\O }{\Omega }
\newcommand{\vp }{\varphi }

\newcommand{\ov}{\overline}
\newcommand{\pa}{\partial}

\newcommand{\wt }{\tilde}

\newcommand{\be}{\begin{equation}}
\newcommand{\ee}{\end{equation}}

\title[symmetry for quasi-linear parabolic
problems]{Asymptotic symmetry for a class of \\ quasi-linear parabolic problems}

\author[L.\ Montoro]{Luigi Montoro$^*$}
\address{Dipartimento di Matematica
\newline\indent
Universit\`a della Calabria
\newline\indent
Ponte Pietro Bucci 31B, I-87036 Arcavacata di Rende, Cosenza, Italy}
\email{montoro@mat.unical.it}

\author[B.\ Sciunzi]{Berardino Sciunzi$^*$}
\address{Dipartimento di Matematica
\newline\indent
Universit\`a della Calabria
\newline\indent
Ponte Pietro Bucci 31B, I-87036 Arcavacata di Rende, Cosenza, Italy}
\email{sciunzi@mat.unical.it}

\thanks{$^*$Dipartimento di Matematica,
Universit\`a della Calabria,
Ponte Pietro Bucci 31B, I-87036 Arcavacata di Rende, Cosenza, Italy,
E-mail: {\em montoro@mat.unical.it}, {\em sciunzi@mat.unical.it}}

\author[M.\ Squassina]{Marco Squassina$^\dagger$}
\address{Dipartimento di Informatica
\newline\indent
Universit\`a degli Studi di Verona
\newline\indent
C\'a Vignal 2, Strada Le Grazie 15, I-37134 Verona, Italy}
\email{marco.squassina@univr.it}

\thanks{$^\dagger$Departimento di Informatica,
Universit\`a di Verona,
C\`a Vignal 2, Strada Le Grazie 15, I-37134 Verona, Italy.
E-mail: {\em marco.squassina@univr.it}}

\thanks{The authors were partially supported by the
Italian PRIN Research Project 2007: {\em Metodi Variazionali e Topologici
nello Studio di Fenomeni non Lineari}}

\begin{document}
	


\begin{abstract}
We study the symmetry properties of the weak positive solutions to a class of quasi-linear
elliptic problems having a variational structure. On this basis, the asymptotic behaviour of
global solutions of the corresponding parabolic equations is also investigated. In particular,
if the domain is a ball, the elements of the $\omega$ limit set are nonnegative radially symmetric
solutions of the stationary problem.
\end{abstract}
\maketitle

\medskip
\begin{center}
\begin{minipage}{11cm}
\footnotesize
\tableofcontents
\end{minipage}
\end{center}

\medskip

\section{Introduction and main results}
Let $\O\subset\R^n$ be a smooth bounded domain and $1<p< \infty$. The goal of this
paper is to study the asymptotic symmetry properties for a class of global solutions
of the following quasi-linear parabolic problem
\begin{equation}
	\label{problema}
	\tag{$E$}
\begin{cases}
u_t-{\rm div}(a(u)|\nabla u|^{p-2}\nabla u)+\frac{a'(u)}{p}|\nabla u|^p=f(u) & \text{in $(0,\infty)\times\Omega$,} \\
u(0,x)=u_0(x)  & \text{in $\Omega$,} \\
u(t,x)=0  & \text{in $(0,\infty)\times\partial\Omega$.}
\end{cases}
\end{equation}
The adoption of the $p$-Laplacian operator inside the diffusion term arises
in various applications where the standard linear heat operator $u_t-\Delta$
is replaced by a nonlinear diffusion with gradient dependent diffusivity. These models
have been used in the theory of non-Newtonian filtration fluids, in turbulent flows in porous media
and in glaciology (cf.~\cite{AE}).
In the following we will assume that $a\in C^2_{{\rm loc}}(\R)$ and there exists a positive constant $\eta$
such that $a(s) \geq \eta>0$ for all $s\in\R^+$ and that $f$ is a locally lipschitz continuous in $[0, \infty)$,
which satisfies some additional positivity conditions. The nontrivial (positive) stationary
solutions of the above problem must be solutions of the following elliptic equation
\begin{equation}\label{eq:pe}\tag{$S$}
\begin{cases}
 -{\rm div}(a(u) |\n u|^{p-2}\n u)+\frac{a'(u)}{p}|\n u|^p = f(u)  & \text {in $\O$}, \\
u>0 & \text{in $\O$}, \\
u=0 &\text{on $\pa \O$}.
\end{cases}
\end{equation}
This class of problems has been intensively studied with respect to existence, nonexistence and multiplicity
via non-smooth critical point theory.
For a quite recent survey paper, we refer the interested reader to~\cite{Sq} and to the references therein.
Already in the investigation of the qualitative properties for the
pure $p$-Laplacian case $a\equiv 1$, one has to face nontrivial difficulties
mainly due to the lack of regularity of the solutions of problem~\eqref{eq:pe}.
As known, the maximal regularity of bounded solutions in the interior of the
domain is $C^{1,\alpha}(\O)$ (see~\cite{Di,T}). Also, since we are assuming the domain to be smooth,
the $C^{1,\alpha}$ regularity assumption up to the boundary follows by~\cite{Li}.
In some sense, the problem is singular (for $1<p<2$) and degenerate (for $p>2$) due to the different behaviour
of the weight $|\nabla u|^{p-2}$.

\begin{definition}
We denote by ${\mathcal S}_{x_1}$ the set of nontrivial weak $C^{1,\a}(\overline{\O})$
solutions $z$ of problem~\eqref{eq:pe} which are symmetric and non-decreasing in the $x_1$-direction\footnote{As customary we  consider the case of a domain which is symmetric with
respect to the hyperplane $\{x_1=0\}$, and we mean that the solution $z$ is non-decreasing in the $x_1$-direction for $x_1<0$. While it is non-increasing for $x_1>0$.}.
We denote by ${\mathcal R}$ the set of nontrivial weak $C^{1,\a}(\overline{\O})$
solutions $z$ of problem~\eqref{eq:pe} which are radially symmetric and radially decreasing.
\end{definition}

\vskip2pt
\noindent
The first result of the paper, regarding the stationary problem, is the following

\begin{theorem}
	\label{main1}
  Assume that $f$ is strictly positive in $(0, \infty )$ and $\O$ is strictly convex with respect
to a direction, say $x_1$, and symmetric with
respect to the hyperplane $\{x_1=0\}$. Then, a weak $C^{1,\a}(\overline{\O})$
solution $u$ of problem~\eqref{eq:pe} belongs to ${\mathcal S}_{x_1}$.
In addition, if $\Omega$ is a ball, then $u$ belongs to ${\mathcal R}$.
\end{theorem}

Following also some ideas in~\cite{DS1}, the main point in proving the
above result is providing in this framework a suitable summability
for the weight $|\nabla u|^{-1}$, allowing to prove that the set of critical points of $u$ has actually
zero Lebesgue measure.

\begin{definition}
	\label{defsolutt}
Given $u_0\in W^{1,p}_0(\Omega)$ with $u_0\geq 0$ a.e.,
we write $u_0\in {\mathcal G}$, if there exists a function
\begin{equation}
	\label{regrequirm}
u\in C([0,\infty);W^{1,p}_0(\Omega,\R^+)),\quad
u_t\in L^2([0,\infty);L^2(\Omega)),\qquad u(0)=u_0,
\end{equation}
solving the problem
\begin{align*}
	\int_0^T\int_{\Omega} u_t\varphi dxdt &+\int_0^T\int_{\Omega}a(u)|\nabla u|^{p-2}\nabla u\cdot\nabla \varphi dx dt \\
&	+\int_0^T\int_{\Omega}\frac{a'(u)}{p}|\nabla u|^{p}\varphi dx dt=\int_0^T\int_{\Omega}f(u)\varphi dx dt,
\quad\forall\varphi\in C^\infty_c(Q_T),
\end{align*}
for any $T>0$, where $Q_T=\Omega\times [0,T]$
and satisfying the energy inequality
\begin{equation}
	\label{energyinequality}
{\mathcal E}(u(t))+\int_s^t\int_{\Omega} |u_t(\tau)|^2dxd\tau \leq {\mathcal E}(u(s)),
\quad\text{for all $t>s\geq 0$,}
\end{equation}
where the energy functional is defined as
$$
{\mathcal E}(u(t))=\frac{1}{p}\int_{\Omega}a(u(t))|\nabla u(t)|^pdx-
\int_{\Omega}F(u(t))dx,\qquad
F(s)=\int_0^s f(\tau)d\tau.
$$
\end{definition}

As we learn from a (classical) work of Tsustumi~\cite[Theorems 1 to 4]{Ts} regarding the pure $p$-Laplacian
case (see also the works~\cite{Is,Zh}), the requirements~\eqref{regrequirm}
in Definition~\ref{defsolutt} are natural. In general,
for the weak solutions of~\eqref{problema} to be globally defined, it is necessary
that the initial datum $u_0$ is chosen sufficiently small. A similar consideration can be done
for the size of the domain $\O$, sufficiently small domains yield global solutions,
while large domains may yield to the appearance of blow-up phenomena.
For well-posedness and H\"older regularity results for quasi-linear parabolic equation, we also
refer the reader to the books~\cite{Di1,Li2}. Finally, concerning the energy inequality~\eqref{energyinequality},
of course smooth solutions of~\eqref{problema} will satisfy the energy identity (namely equality
in~\eqref{energyinequality} in place of the inequality). It is sufficient to multiply~\eqref{problema}
by $u_t$ and, then, integrate in space and time. On the other hand~\eqref{energyinequality} is enough for our purposes
and it seems implicitly automatically satisfied by the Galerkin method yielding the existence and regularity
of solutions, see e.g.~\cite[identity (3.8) and related weak convergences (3.9)-(3.13)]{Ts}.
\vskip2pt
\noindent
The second result of the paper is the following

\begin{theorem}
	\label{main-mid}
Assume that there exists a positive constant $\rho$ such that
\begin{equation}
	\label{segnocondiz}
a'(s)s\geq 0,\quad\text{for all $s\in\R$ with $|s|\geq \rho$},
\end{equation}
and that there exist two positive constants
$C_1,C_2$ and $\sigma\in [1,p^*-1)$ with $p>\frac{2n}{n+2}$, such that
\begin{equation}
	\label{growth-f}
|f(s)|\leq C_1+C_2|s|^\sigma,\qquad\text{for all $s\in\R$}
\end{equation}
Then, the following facts hold.
\begin{itemize}
\item[${\rm (a)}$]
Assume that $f$ is strictly positive in $(0, \infty )$ and $\O$ is strictly convex with respect
to a direction, say $x_1$, and symmetric with respect to the hyperplane $\{x_1=0\}$.
Let $u_0\in {\mathcal G}$ and let $u:[0,\infty)\times\Omega\to\R^+$ be
the corresponding solution of~\eqref{problema}. Then, for any diverging
sequence $(\tau_j)\subset\R^+$ there exists a diverging
sequence $(t_j)\subset\R^+$ with $t_j\in[\tau_j,\tau_j+1]$ such that
$$
u(t_j)\to z\quad\text{strongly in $W^{1,p}_0(\O)$ as $j\to\infty$},
$$
where either $z=0$ or $z\in{\mathcal S}_{x_1}$ (if $\Omega=B(0,R)$ with $R>0$, then
either $z=0$ or $z\in{\mathcal R}$) provided that $z\in L^\infty(\O)$. In addition, for all $\mu_0>0$,
	\begin{equation}
		\label{additionalinfo}
\sup_{\mu\in [0,\mu_0]}\|u(t_j+\mu)-z\|_{L^q(\O)}\to 0\quad\text{as $j\to\infty$},
\end{equation}
for any $q\in [1,p^*)$.
\vskip6pt
\item[${\rm (b)}$] Let $R>0$ and assume that $f\in C^1([0,\infty))$ with $f(0)=0$ and
\begin{equation}
	\label{f-plus-1}
	0<(p-1)f(s)<sf'(s),\qquad\text{for all $s>0$}.
\end{equation}
Furthermore, assume that
\begin{equation}
	\label{f-plus-2}
H'(s)\leq 0\quad\text{for $s>0$},\quad H(s)=(n-p)s-np\frac{\int_0^s f(\tau)d\tau}{f(s)}
,\quad\text{$H(0)=0$.}
\end{equation}
Let $u_0\in {\mathcal G}$ and let $u:[0,\infty)\times B(0,R)\to\R^+$ be
the corresponding solution of
\begin{equation}
	\label{plapl-ev}
\begin{cases}
u_t-\Delta_p u=f(u) & \text{in $(0,\infty)\times B(0,R)$,} \\
u(0,x)=u_0(x)  & \text{in $B(0,R)$,} \\
u(t,x)=0  & \text{in $(0,\infty)\times\partial B(0,R)$.}
\end{cases}
\end{equation}
Then, for any diverging
sequence $(\tau_j)\subset\R^+$ there exists a diverging
sequence $(t_j)\subset\R^+$ with $t_j\in[\tau_j,\tau_j+1]$ such that
$$
u(t_j)\to z\quad\text{strongly in $W^{1,p}_0(\O)$ as $j\to\infty$},
$$
where either $z=0$ or $z$ is the unique positive solution to the problem
\begin{equation}
	\label{plapl-st}
\begin{cases}
-\Delta_p u=f(u) & \text{in $B(0,R)$,} \\
u>0  & \text{in $B(0,R)$,} \\
u=0  & \text{on $\partial B(0,R)$.}
\end{cases}
\end{equation}
In addition, the limit~\eqref{additionalinfo} holds.
\end{itemize}
\end{theorem}

\vskip2pt
\noindent
\begin{remark}
	The sign condition~\eqref{segnocondiz} is often assumed in the current literature on problem~\eqref{eq:pe}
	(and in more general frameworks as well) in dealing with both existence and nonexistence results (see e.g.~\cite{CD,Sq,BBM}).
	We point out that it is, in general, necessary for the mere $W^{1,p}_0(\O)$ solutions to~\eqref{eq:pe} to be bounded in
	$L^\infty(\O)$ (see~\cite{Fr}).
\end{remark}

\vskip2pt
\noindent
Next, we consider a class of initial data corresponding to global solutions which enjoy
some compactness over, say, the time interval $\{t>1\}$.

\begin{definition}
We write $u_0\in {\mathcal A}$ if $u_0\in {\mathcal G}$ and furthermore, the set
$$
K=\big\{u(t): t>1\big\},
$$
is relatively compact in $W^{1,p}_0(\O)$.
For any initial datum $u_0\in W^{1,p}_0(\Omega)$, the $\omega$-limit set of $u_0$ is defined as
\begin{equation*}
	\omega(u_0)=\big\{ z\in W^{1,p}_0(\Omega):
	\text{there is $(t_j)\subset\R^+$ with
	$u(t_j)\to z$ in $W^{1,p}_0(\Omega)$}\big\},
\end{equation*}
where $u(t)$ is the solution of~\eqref{problema} corresponding to $u_0$.
\end{definition}

\vskip4pt
\noindent
The third, and last, result of the paper is the following

\begin{theorem}
	\label{main2}
Assume that $f$ is strictly positive in $(0, \infty )$ with the growth~\eqref{growth-f} and $\O$ is strictly convex with respect
to a direction, say $x_1$, and symmetric with respect to the hyperplane $\{x_1=0\}$.
Then, the following facts hold.
\begin{itemize}
\item[${\rm (a)}$] For all $u_0\in {\mathcal A}$, we have
$$
\omega(u_0)\cap L^\infty(\O)\subset {\mathcal S}_{x_1}.
$$
In particular, the $L^\infty$-bounded elements of the $\omega$-limit set to~\eqref{problema} with $\Omega=B(0,R)$
are zero or  radially symmetric and decreasing
solutions of problem~\eqref{eq:pe}.
\item[${\rm (b)}$]
Assume that $f\in C^1([0,\infty))$ with $f(0)=0$ satisfies
assumptions~\eqref{f-plus-1} and~\eqref{f-plus-2}.
Then, for all $u_0\in {\mathcal A}$, the $\omega$-limit set of problem~\eqref{plapl-ev}
consists of either $0$ or the unique positive solution to the problem~\eqref{plapl-st}.
\end{itemize}
\end{theorem}
\begin{remark}
Quite often, even in the fully nonlinear parabolic case, global
solutions which are uniformly
bounded in $L^\infty$ are considered (see e.g.~\cite[Section 3.1]
{Po}). In these cases,
in our framework, the elements
of the $\omega$-limit set are automatically bounded and, in turn,
belong to $C^{1,\a}(\overline{\O})$. Concerning the $L^
\infty$-global boundedness issue for a class of
degenerate operators, such as the $p$-Laplacian case, we refer the
reader to the work of Lieberman~\cite{Li1}, in particular~
\cite[Theorem 2.4]{Li1}, where he proves that
$$
\sup_{(t,x)\in [0,\infty)\times\O}|u(t,x)|<\infty,
$$
provided that suitable growth conditions hold on the parabolic operator as well
as on the nonlinearity,
which satisfy a typical super-linearity condition, reading as
$$
f(s)s\geq (a_0+\alpha)F(s)-c_1,\quad F(s)\geq s^{2+\alpha}-c_0,\quad s
\in\R,
$$
for suitable positive constants $a_0,c_0,c_1$ and $\alpha$.
\end{remark}

\begin{remark}
Assume that $\Omega$ is a star-shaped domain and consider the problem with the critical power nonlinearity	
\begin{equation}
	\label{crit-eq}
\begin{cases}
 -{\rm div}(a(u) |\n u|^{p-2}\n u)+\frac{a'(u)}{p}|\n u|^p = u^{p^*-1}  & \text {in $\O$}, \\
u>0 & \text{in $\O$}, \\
u=0 &\text{on $\pa \O$}.
\end{cases}
\end{equation}	
Assuming the sign condition
$$
a'(s)\geq 0,\quad\text{for all $s\geq 0$},
$$ 	
it is known that problem~\eqref{crit-eq} does not admit any solution (cf.~\cite{PS,DMS}).
In turn, any uniformly bounded global solution to the problem	
\begin{equation*}
\begin{cases}
u_t-{\rm div}(a(u)|\nabla u|^{p-2}\nabla u)+\frac{a'(u)}{p}|\nabla u|^p=u^{p^*-1}  & \text{in $(0,\infty)\times\Omega$,} \\
u(0,x)=u_0(x)  & \text{in $\Omega$,} \\
u(t,x)=0  & \text{in $(0,\infty)\times\partial\Omega$}
\end{cases}
\end{equation*}
must vanish along diverging sequences
$(t_j)\subset\R^+$, $u(t_j)\to 0$ in $W^{1,p}_0(\O)$ as $j\to\infty$.
\end{remark}

\begin{remark}
Theorems~\ref{main1},~\ref{main-mid},~\ref{main2} are new already
in the non-degenerate case $p=2$ since of the presence of the coefficient $a(\cdot)$, in which case the solutions are expected
to be very regular for $t>0$.
\end{remark}
\vskip2pt
\noindent

We do not investigate here conditions under which one can characterize a class
of initial data which guarantee global solvability with the additional information
of compactness of the trajectory into $W^{1,p}_0(\O)$.
In the semi-linear case $p=2$ with a power type nonlinearity $f(u)=|u|^{m-1}u$, $m>1$,
we refer to~\cite{CL,Qu,Qu1} for apriori estimates and smoothing properties in $C^1(\Omega)$
of the solutions for positive times. About the convergence to
nontrivial solutions to the stationary problem along some suitable diverging
time sequence $(t_j)\subset\R^+$, we also refer to~\cite{GW}
for a detailed analysis of the sets of initial data $u_0\in H^1_0(\O)$ yielding
to vanishing and non-vanishing global solutions as well as initial data for which
the solutions blow-up in finite time.
In particular it is proved that the stabilization
towards nontrivial equilibria is a borderline case, in the sense that the set of initial data
corresponding to non-vanishing global solution is precisely the boundary of the (closed) set of data
yielding global solutions.
In conclusion, in general, at least four different type of behaviour may occur in these problems:
blow up in finite time, global vanishing solution, global non-vanishing solution (converging to equilibria) and
finally global solution blowing up in infinite time (see also~\cite{NST}).
In our general framework, also due to the degenerate nature of the problem,
this classification seems quite hard to prove, so we focus on the third case.
In the $p$-Laplacian case $a\equiv 1$, we refer
the reader to~\cite{Li1} for the study of apriori estimates and convergence
to equilibria for global solutions. Our approach is based on the independent study
of the symmetry properties of positive stationary solutions via a suitable weak comparison principle
allowing to apply the Alexandrov-Serrin moving plane technique
in symmetric domains (see also~\cite{DP,DS1,DS2} for similar results
in the case $a=1$). Then, since the problem clearly admits a variational structure and the
energy functional ${\mathcal E}:W^{1,p}_0(\Omega)\to\R$ defined by
$$
{\mathcal E}(u(t))=\frac{1}{p}\int_{\Omega}a(u(t))|\nabla u(t)|^pdx-\int_{\Omega}F(u(t))dx,\quad t>0,
\quad F(s)=\int_0^sf(\tau)d\tau,
$$
is decreasing along a smooth solution $u(t)$, the global solutions have to approach
stationary states along suitable diverging sequences $(t_j)\subset\R^+$. In pursuing this target
we also make use of some nontrivial compactness result proved in~\cite{CD} in the study of the
stationary problem. It is known that, in general,
it is not possible to get the convergence result along the whole trajectory, namely as $t\to\infty$ (see~\cite{PoSi})
unless the nonlinearity $f$ is an analytic function (see~\cite{Je}).

For a general survey paper on the asymptotic
symmetry of the solutions to general (not just those with a Lyapunov functional)
nonlinear parabolic problems,
we refer to the recent work of P.~Pol{\'a}{\v{c}}ik~\cite{Po}
where various different approaches to the study of the problem are discussed.
\bigskip
\vskip8pt
\begin{center}\textbf{Plan of the paper.}\end{center}
\vskip2pt
\noindent
In Section~\ref{stationary} we study the regularity properties of the weak positive
solutions to~\eqref{eq:pe}. In Section~\ref{parabolic} we obtain some properties related to the asymptotic behaviour of solutions to the parabolic problem~\eqref{problema}. Finally, in Section~\ref{prove} we complete the proof of the main results of the paper.
\bigskip
\vskip8pt
\begin{center}\textbf{Notations.}\end{center}
\begin{enumerate}
\item For $n\geq 1$, we denote by $|\cdot|$ the euclidean norm in $\R^n$.
\item $\R^+$ (resp.\ $\R^-$) is the set of positive (resp.\ negative) real values.
\item For $p>1$ we denote by $L^p(\R^n)$ the space of measurable functions $u$ such that
$\int_{\O}|u|^pdx<\infty$. The norm $(\int_{\O}|u|^pdx)^{1/p}$ in $L^p(\O)$ is denoted by $\|\cdot\|_{L^p(\O)}$.
\item For $s\in\N$, we denote by $H^s(\O)$ the Sobolev space of functions $u$ in $L^2(\O)$
having generalized partial derivatives $\partial_i^k u$ in $L^2(\O)$ for all $i=1,\dots, n$ and
any $0\leq k \leq s$.
\item The norm $(\int_{\O}|u|^pdx+ \int_{\O}|\nabla u|^pdx)^{1/2}$ in $W^{1,p}_0(\O)$ is denoted by $\|\cdot\|_{W^{1,p}_0(\O)}$.
\item We denote by $C_0^{\infty}(\O)$ the set of smooth compactly supported functions in $\O$.
\item We denote by $B(x_0,R)$ a ball of center $x_0$ and radius $R$.
\item We denote $D^2u$ the Hessian matrix of $u$ and  $| D^2 u|^2 \equiv \sum_{i=1}^n |\n u_i|^2$.
\item We denote by $\mathcal{L}(E)$ the Lebesgue measure of the set $E\subset\R^n$.
\end{enumerate}
\medskip

\section{Symmetry for stationary solutions}
\label{stationary}

We consider weak $C^{1,\a}(\overline{\O})$ solutions
to~\eqref{eq:pe}. We recall that we shall assume that
\begin{itemize}\label{eq:hp}
  \item [$({\rm i})$] $f$ is locally lipschitz continuous in $[0, \infty)$;
  \item [$({\rm ii})$] For any given $\t>0$, there exists a positive
  constant $K$ such that $f(s)+Ks^q\geq0$ for some $q \geq p-1$ and
  for any $s \in [0, \t]$. Observe that this implies $f(0)\geq0$;
  \item [$({\rm iii})$] $a\in C^2_{{\rm loc}}(\R)$ and there exists $\eta>0$
such that $a(t) \geq \eta>0$;
\end{itemize}

As pointed out in the introduction, if we assume that the solution is bounded,
the $C^{1,\alpha}$ regularity up to the boundary follows by~\cite{Di,T,Li}.
Also hypothesis (iii) ensures the applicability of the Hopf boundary lemma (see~\cite{PSB,PSZ}).

\subsection{Gradients summability}
In weak form, our problem reads as
\begin{equation}
	\label{eqdebole}
\int_\O a(u) |\n u|^{p-2}\n u\cdot\n \vp dx+\frac 1p\int_\O
a'(u)|\n u|^p\vp dx= \int_\O f(u)\vp dx,\quad\forall \vp\in C^\infty_c(\Omega).
\end{equation}
Define, as usual, the critical set $Z_u$ of $u$ by setting
\begin{equation}\label{eq:Z_u}
Z_u=\big \{x \in \O : \n u(x)=0 \big \}
\end{equation}
Note that the importance of critical set $Z_u$ is due to the
fact that it is exactly the set where our operator is degenerate.
By Hopf Lemma (cf.\ \cite{PSB,PSZ}), it follows that
\begin{equation}
	\label{eq:Z_u and boundary}
Z_u\cap \pa\O= \emptyset.
\end{equation}
We want to point out that, by standard regularity results,  $u\in
C^2_{{\rm loc}}(\O\setminus Z_u)$. For functions $\vp\in C^\infty_c(\O\setminus Z_u)$,
let us consider the test function  $\vp_i=\partial_{x_i}\vp$ and denote also
$u_i= \partial_{x_i} u$, for all $i=1,\dots,n$.
With this choice in~\eqref{eqdebole}, integrating by part, we get
\begin{align}\label{eq:linearizzato}\nonumber
\int_\O a(u) |\n u|^{p-2}(\n u_i, \n \vp)&+(p-2)\int_\O a(u)|\n u|^{p-4}(\n
u, \n u_i)(\n u, \n \vp)dx  \notag \\
&+ \int_\O  a'(u)|\n u|^{p-2}(\n u, \n
\vp) u_i dx    \\
&+\int_\O \frac 1p a''(u)|\n u|^p u_i \vp + \int_\O a'(u)|\n
u|^{p-2}(\n u, \n u_i) \vp \notag \\
&-\int_\O f'(u)u_i \vp=0,   \notag
\end{align}
that is, in such a way, we have defined the linearized operator
$L_u(u_i, \vp)$ at a fixed solution $u$ of \eqref{eq:pe}. Then we
can write   equation \eqref{eq:linearizzato} as
\begin{equation}\label{L_u}
L_u(u_i,\vp)=0, \qquad \forall \vp  \in C^\infty_c(\O\setminus Z_u).
\end{equation}

\vskip2pt
\noindent
In the following, we repeatedly use Young's inequality in this form
$$
ab\leq  \delta a^2+C(\d) b^2\qquad\text{for all $a,b\in\R$ and $\d>0$}.
$$
\vskip3pt
We can now state and prove the following

\begin{proposition}
	\label{pr:Summability}
	Let $u \in C^{1, \a}(\ov{\O})$ be a solution to problem \eqref{eq:pe}. Assume that $f$ is locally
lipschitz continuous, $a \in C^2_{{\rm loc}}(\R)$ and there exists a positive constant $\eta$
such that $a(s) \geq \eta>0$ for all $s\in\R^+$.
 Assume that $\O$ is a bounded and smooth domain of $\R^n$.
Then
\begin{equation}\label{eq:Tesi 1}
\int_{\O\setminus \{u_i=0\}} \frac{|\n u|^{p-2}}{|y-x|^\g}\frac{|\n
u_i|^2}{|u_i|^\b} \,dx\leq\mathcal{C},
\end{equation}
where $0 \leq \b < 1$, $\g < n-2$ $(\g=0$ if $n=2)$, $1<p<\infty$
and the positive constant $\mathcal{C}$ does not depend on $y$. In particular,
we have
\begin{equation}\label{eq:Tesi 2}
\int_{\O\setminus \{\n u=0\}} \frac{|\n
u|^{p-2-\b}||D^2u||^2}{|y-x|^\g} \,dx\leq\mathcal{\tilde{C}},
\end{equation}
for a positive constant $\mathcal{\tilde{C}}$ not depending on $y$.
\end{proposition}

\begin{proof}
For all $\e>0$, let us define the piecewise smooth function $G_\eps:\R\to\R$ by setting
\begin{equation}\label{eq:G}
G_\e(t)=\begin{cases} t & \text{if  $|t|\geq 2\e$}, \\
2t-2\e& \text{if $\e\leq t\leq2\e$}, \\
2t+2\e& \text{if  $-2\e\leq t\leq-\e$},
\\ 0 & \text{if $|t| \leq \e$}.
\end{cases}
\end{equation}
Let us choose $E \subset \subset \O$ and a positive function $\psi \in C^\infty_c(\O)$,
such that the support of $\psi$ is compactly contained in $\O$,
$\psi \geq 0$ in $\O$ and $\psi \equiv 1$ in $E$. Let us set
\begin{equation}
	\label{eq:vp}
\vp_{\e,y}(x)=\frac{G_\e(u_i(x))}{|u_i(x)|^\b}
\frac{\psi(x)}{|y-x|^\g}
\end{equation}
where $0 \leq \b <1$, $\g < n-2$ ($\gamma=0$ for $n=2$).  Since $\vp_{\e,y} $ vanishes in a
neighborhood of each critical point, it follows that $\vp_{\e,y} \in
C^2_c(\O \setminus Z_u)$ and hence we can use it as a test function
in \eqref{eq:linearizzato}, getting the following result
\begin{align*}
&  \int_\O  \frac{a(u)}{|y-x|^\g}\frac{|\n u|^{p-2}}{|u_i|^\b}\Big
(G'_\e(u_i)-\b \frac{G_\e(u_i)}{u_i} \Big)\psi |\n u_i|^2dx\\
 &+\int_\O (p-2)\frac{a(u)}{|y-x|^\g}\frac{|\n
u|^{p-4}}{|u_i|^\b}\Big (G'_\e(u_i)-\b \frac{G_\e(u_i)}{u_i}
\Big)\psi(\n u, \n u_i)^2 dx\\ \nonumber & +\int_\O
 \frac{a'(u)}{|y-x|^\g}\frac{|\n u|^{p-2}}{|u_i|^\b}\Big
(G'_\e(u_i)-\b \frac{G_\e(u_i)}{u_i} \Big)\psi u_i (\n u, \n
u_i)dx  \\
&+ \int_\O \frac{a(u)}{|y-x|^\g}|\n
u|^{p-2}\frac{G_\e(u_i)}{|u_i|^\b}(\n u_i,\n \psi)dx\\
 & +\int_\O(p-2)\frac{a(u)}{|y-x|^\g}|\n
u|^{p-4}\frac{G_\e(u_i)}{|u_i|^\b}(\n u, \n u_i)(\n u, \n \psi)dx\\
 & + \int_\O \frac{a'(u)}{|y-x|^\g}|\n u|^{p-2}u_i
\frac{G_\e(u_i)}{|u_i|^\b}(\n u, \n \psi)dx \\
& + \int_\O
a(u)|\n u|^{p-2}\frac{G_\e(u_i)}{|u_i|^\b}\psi(\n u_i,
\n_x(\frac{1}{|y-x|^\g}))dx\\  & +\int_\O (p-2)a(u)|\n
u|^{p-4}\frac{G_\e(u_i)}{|u_i|^\b}\psi(\n u, \n u_i)(\n
u,\n_x(\frac{1}{|y-x|^\g}))dx\\ & + \int_\O  a'(u)|\n
u|^{p-2}u_i\frac{G_\e(u_i)}{|u_i|^\b}\psi(\n
u,\n_x(\frac{1}{|y-x|^\g}))dx   \\
& +\int_\O\frac 1p
a''(u)|\n
u|^p u_i \frac{G_\e(u_i)}{|u_i|^\b}\frac{\psi}{|y-x|^\g}dx\\
& +\int_\O  a'(u)|\n u|^{p-2}(\n u, \n
u_i)\frac{G_\e(u_i)}{|u_i|^\b}\frac{\psi}{|y-x|^\g}dx  = \int_\O
f'(u)u_i\frac{G_\e(u_i)}{|u_i|^\b}\frac{\psi}{|y-x|^\g}dx
\end{align*}
Let us denote each term of the previous equation in
a useful way for the sequel, that is
\begin{align}\label{eq:A_1}
& A_1=\int_\O  \frac{a(u)}{|y-x|^\g}\frac{|\n
u|^{p-2}}{|u_i|^\b}\Big (G'_\e(u_i)-\b \frac{G_\e(u_i)}{u_i}
\Big)\psi |\n u_i|^2dx;\\ \nonumber
&A_2=\int_\O
(p-2)\frac{a(u)}{|y-x|^\g}\frac{|\n u|^{p-4}}{|u_i|^\b}
\Big(G'_\e(u_i)-\b \frac{G_\e(u_i)}{u_i} \Big)\psi(\n u, \n u_i)^2 dx;
\\ \nonumber
\nonumber & A_3=\int_\O \frac{a'(u)}{|y-x|^\g}\frac{|\n
u|^{p-2}}{|u_i|^\b}\Big (G'_\e(u_i)-\b \frac{G_\e(u_i)}{u_i}
\Big)\psi u_i (\n u, \n
u_i)dx;\\
\nonumber & A_4= \int_\O \frac{a(u)}{|y-x|^\g}|\n
u|^{p-2}\frac{G_\e(u_i)}{|u_i|^\b}(\n u_i,\n \psi)dx;\\
\nonumber& A_5=\int_\O(p-2)\frac{a(u)}{|y-x|^\g}|\n
u|^{p-4}\frac{G_\e(u_i)}{|u_i|^\b}(\n u, \n u_i)(\n u, \n \psi)dx;\\
\nonumber & A_6= \int_\O \frac{a'(u)}{|y-x|^\g}|\n u|^{p-2}u_i
\frac{G_\e(u_i)}{|u_i|^\b}(\n u, \n \psi)dx;
\\ \nonumber& A_7= \int_\O
a(u)|\n u|^{p-2}\frac{G_\e(u_i)}{|u_i|^\b}\psi(\n u_i,
\n_x(\frac{1}{|y-x|^\g}))dx;\\ \nonumber & A_8=\int_\O (p-2)a(u)|\n
u|^{p-4}\frac{G_\e(u_i)}{|u_i|^\b}\psi(\n u, \n u_i)(\n
u,\n_x(\frac{1}{|y-x|^\g}))dx;
\\ \nonumber& A_9= \int_\O a'(u)|\n
u|^{p-2}u_i\frac{G_\e(u_i)}{|u_i|^\b}\psi(\n
u,\n_x(\frac{1}{|y-x|^\g}))dx;
\\ \nonumber &A_{10}=\int_\O \frac 1p a''(u)|\n u|^p
u_i \frac{G_\e(u_i)}{|u_i|^\b}\frac{\psi}{|y-x|^\g}dx;
\\ \nonumber&
A_{11}=\int_\O  a'(u)|\n u|^{p-2}(\n u, \n
u_i)\frac{G_\e(u_i)}{|u_i|^\b}\frac{\psi}{|y-x|^\g}dx;
\\ \nonumber &
N= \int_\O f'(u)u_i\frac{G_\e(u_i)}{|u_i|^\b}\frac{\psi}{|y-x|^\g}dx.
\end{align}
Then we have rearranged the equation as
\begin{equation}\label{eq:linearizzato2}
\sum_{i=1}^{11}A_i=N
\end{equation}
Notice that, since $0\leq\b <1$, for all $t\in\R$ and $\eps>0$ we have
\begin{equation*}
 G'_\e(t)- \frac{\b G_\e(t)}{t} \geq 0,\qquad
  \lim_{\eps\to 0}\Big(G'_\e(t)- \frac{\b G_\e(t)}{t}\Big)=1-\b.
\end{equation*}
From now on, we will denote
$$
\tilde{G}_\e(t)=  G'_\e(t)- \b\frac{G_\e(t)}{t},\qquad \text{for all $t\in\R$ and $\eps>0$}.
$$
From equation~\eqref{eq:linearizzato2} one has
$$
A_1+A_2\leq \sum_{i=3}^{11}|A_i|+|N|.
$$
We shall distinguish the proof into two cases.
\vskip2pt
\noindent
{\bf Case I:} $\mathbf{p \geq 2}$. This trivially implies  $A_2\geq 0$, and hence
\begin{equation}
	\label{eq:case i}
	A_1\leq A_1+A_2\leq \sum_{i=3}^{11}|A_i|+|N|.
	\end{equation}
\vskip2pt
\noindent
{\bf Case II:} $\mathbf{1<p<2}$.
By Schwarz inequality, of course, it follows
$$
|\n u|^{p-4}(\n u, \n u_i)^2 \leq |\n u|^{p-2}|\n u_i|^2.
$$
In turn, since $1<p<2$, this implies
$$
(p-2)a(u)\frac{\tilde G_\e(u_i)}{|u_i|^\b}
\frac{\psi |\n u|^{p-4}(\n u,\n u_i)^2}{|y-x|^\g}\geq(p-2)a(u)\frac{\tilde G_\e(u_i)}{|u_i|^\b}
\frac{\psi |\n u|^{p-2}|\n u_i|^2}{|y-x|^\g},
$$
so that $(p-2)A_1\leq A_2$, yielding
\begin{equation}\label{eq:case ii}
A_1 \leq \frac{1}{p-1}(A_1+A_2)\leq \frac{1}{p-1}\sum_{i=3}^{11}|A_i|+|N|.
\end{equation}
In both cases, in view of~\eqref{eq:case i} and~\eqref{eq:case ii},
we want to estimates the terms in the sum
\begin{equation}
	\label{eq:esitmates}
\sum_{i=3}^{11}|A_i|+|N|.
\end{equation}
Let us start by estimating the terms $A_i$ in the sum~\eqref{eq:esitmates}.
Concerning $A_3$, we have
\begin{align*}
|A_3|&\leq \int_\O
 \frac{|a'(u)|}{|y-x|^\g}\frac{|\n u|^{p-2}}{|u_i|^\b}
\tilde{G}_\e(u_i)\psi |u_i| |\n u| |\n u_i|dx\\
&\leq C_3 \int_\O \frac{1}{|y-x|^\g}\frac{|\n u|^{p-1}}{|u_i|^\b}
\tilde{G}_\e(u_i)\psi |u_i|  |\n u_i|dx\\
&\leq C_3 \left[ \d \int_\O \frac{|\n
u|^{p-2}}{|y-x|^\g}\frac{\tilde{G}_\e(u_i)}{|u_i|^\b}
\psi |\n u_i|^2dx+C_\d\int_\O\frac{|\n u|^{p-1}}{|y-x|^\g} \psi\frac{\tilde{G}_\e(u_i)}{|u_i|^{\b-2}}dx \right]\\
&\leq \frac{C_3 \d}{\eta} A_1+K_{3}(\d),
\end{align*}
where we used that
$$
|\n u|^{p-1}\psi\frac{\tilde{G}_\e(u_i)}{|u_i|^{\b-2}}\leq C,
$$
where $C$ is a positive constant independent of $\e$ and
$C_3$ is a positive constant independent of $y$. Moreover recall that $0 \leq \b <1$ and that $u \in C^{1, \a}(\ov{\O})$. Also
\begin{equation*}
|A_4|\leq \int_\O \frac{a(u)}{|y-x|^\g}|\n
u|^{p-2}\frac{|G_\e(u_i)|}{|u_i|^\b}|\n u_i||\n \psi|dx\leq C_4,
\end{equation*}
where
$$
{\displaystyle \frac{1}{|y-x|^\g}\frac{|\n
u|^{p-2}}{|u_i|^{\b-1}}\frac{|G_\e(u_i)|}{|u_i|}|\n u_i||\n \psi|}\in L^\infty(\Omega),
$$
since $|\n u_i|$ is bounded in a neighborhood
of the boundary by Hopf Lemma,  $\g-2 < n$, $0 \leq \b <1$ and the constant $C_4$
is independent of $y$. For the same reasons, we also have
\begin{align*}
|A_5|&\leq \int_\O\frac{a(u)}{|y-x|^\g}|\n
u|^{p-2}\frac{|G_\e(u_i)|}{|u_i|^\b}|\n u_i||\n \psi|dx\leq  C_5, \\
|A_6|&\leq \int_\O  \frac{|a'(u)|}{|y-x|^\g}|\n u|^{p-1}
\frac{|G_\e(u_i)|}{|u_i|^{\b-1}} |\n \psi|dx\leq C_6,
\end{align*}
for some positive constants $C_5$ and $C_6$ independent of $y$. Furthermore,
for a positive constant $C_7$ independent of $y$, we have
\begin{align*}
|A_7|&\leq \int_\O
a(u)|\n u|^{p-2}\frac{|G_\e(u_i)|}{|u_i|^\b}\psi|\n u_i|
\big|\n_x \frac{1}{|y-x|^\g}\big|dx\\
&\leq C_7\int_\O a(u)|\n
u|^{p-2}\frac{|G_\e(u_i)|}{|u_i|^\b}\psi|\n u_i|
\frac{1}{|y-x|^{\g+1}}dx\\
&\leq C_7\d\int_\O\frac{a(u)}{|y-x|^{\g}}\frac{|\n
u|^{p-2}}{|u_i|^\b} \psi \frac{|G_\e(u_i)|}{|u_i|}|\n
u_i|^2dx \\
&+C(\d)\int_\O a(u)|\n u|^{p-1}\frac{|G_\e(u_i)|}{|u_i|}\frac{1}{|y-x|^{\g+2}}dx\\
&\leq C_7\d\int_\O\frac{a(u)}{|y-x|^{\g}}\frac{|\n
u|^{p-2}}{|u_i|^\b} \psi \frac{|G_\e(u_i)|}{|u_i|}|\n u_i|^2dx+K_7(\d)
\end{align*}
where we used Young's inequality, $\g -2<n$ and $0 \leq \b <1$.
In a similar fashion,
\begin{align*}
|A_8|&\leq \int_\O |p-2|a(u)|\n
u|^{p-2}\frac{|G_\e(u_i)|}{|u_i|^\b}\psi |\n u_i|\big|\n_x\frac{1}{|y-x|^\g}\big|dx\\
&\leq C_8\d\int_\O\frac{a(u)}{|y-x|^{\g}}\frac{|\n
u|^{p-2}}{|u_i|^\b} \psi \frac{G_\e(u_i)}{u_i}|\n u_i|^2dx+K_8(\d)
\end{align*}
as well as
\begin{equation*}
|A_9| \leq \int_\O  |a'(u)||\n
u|^{p-1}\frac{|G_\e(u_i)|}{|u_i|^{\b-1}}\psi\big|\n_x \frac{1}{|y-x|^\g}\big|dx \leq C_9.
\end{equation*}
for some positive constants $C_8, C_9$ independent of $y$.
We get an upper bound for the last terms as well
\begin{equation*}
|A_{10}| \leq \frac 1p\int_\O |a''(u)||\n u|^p
 \frac{|G_\e(u_i)|}{|u_i|^{\b-1}}\frac{\psi}{|y-x|^\g}dx \leq C_{10},
\end{equation*}
with $C_{10}$ independent
of $y$ and where we have also used the fact that $a \in
C^{2}_{{\rm loc}}(\R)$. In the same way, it holds
\begin{align*}
	|A_{11}|&\leq \int_\O  |a'(u)||\n
u|^{p-1}\frac{|G_\e(u_i)|}{|u_i|^\b}|\n
u_i|\frac{\psi}{|y-x|^\g}dx\\
&\leq C_{11}
\d\int_\O\frac{1}{|y-x|^\g}\frac{|\n
u|^{p-2}}{|u_i|^\b}\frac{G_\e(u_i)}{u_i}\psi|\n
u_i|^2dx+C(\d)\int_\O\frac{|\n
u|^p}{|y-x|^\g}\frac{\psi}{|u_i|^{\b-1}}\\
&\leq
\frac{C_{11}\d}{\eta} \int_\O\frac{a(u)}{|y-x|^\g}\frac{|\n
u|^{p-2}}{|u_i|^\b}\frac{G_\e(u_i)}{u_i}\psi|\n u_i|^2dx+K_{11}(\d)
\end{align*}
and
\begin{align*}
|N| \leq \int_\O |f'(u)|\frac{|G_\e(u_i)|}{|u_i|^{\b-1}}\frac{\psi}{|y-x|^\g}dx \leq C_N,
\end{align*}
where the last inequality holds true since $f$ is locally lipschitz
continuous and where $C_{11}$ and $C_N$ are constants independent
of $y$. Then, by these estimates above and by equations \eqref{eq:case i},
\eqref{eq:case ii} and \eqref{eq:esitmates} we write
\begin{equation}\label{eq:esitmates1}
A_1\leq \mathcal{D}\sum_{i=3}^{11}|A_i|+|N|\leq
\mathcal{S}\d A_1
+\mathcal{M}\d\int_\O\frac{a(u)}{|y-x|^{\g}}\frac{|\n
u|^{p-2}}{|u_i|^\b} \psi \frac{G_\e(u_i)}{u_i}|\n u_i|^2dx +
\mathcal{C}_\d,
\end{equation}
where we have set
\begin{align*}
 \mathcal{D}&= \max\Big\{1, \frac{1}{p-1}\Big\},
\quad \mathcal{S}=\mathcal{D} \frac{C_3}{\eta},
\quad\mathcal{M}=\mathcal{D}\max\Big\{C_7 , C_8 ,
\frac{C_{11}}{\eta}\Big\}\\
 \mathcal{C}_\d&=\max\big\{K_{3}(\d),K_{7}(\d),
K_{8}(\d), K_{11}(\d), C_4,C_5,C_6,C_9,C_N\big\}.
\end{align*}
Then from equations \eqref{eq:A_1} and \eqref{eq:esitmates1} one has
\begin{align*}\nonumber
&(1-\mathcal{S} \d)\int_\O  \frac{a(u)}{|y-x|^\g}\frac{|\n
u|^{p-2}}{|u_i|^\b}\left (G'_\e(u_i)-\b \frac{G_\e(u_i)}{u_i}
\right)\psi |\n u_i|^2dx\\
&\leq\mathcal{M}\d\int_\O\frac{a(u)}{|y-x|^{\g}}\frac{|\n
u|^{p-2}}{|u_i|^\b} \psi \frac{G_\e(u_i)}{u_i}|\n u_i|^2dx +
\mathcal{C}_\d,
\end{align*}
namely
\begin{equation}
(1-\mathcal{S} \d)\int_\O  \frac{a(u)}{|y-x|^\g}\frac{|\n
u|^{p-2}}{|u_i|^\b}\left
[G'_\e(u_i)-\left(\b+\frac{\mathcal{M}\d}{(1-\mathcal{S} \d)}\right)
\frac{G_\e(u_i)}{u_i} \right]\psi |\n u_i|^2dx\leq\mathcal{C}_\d
\end{equation}
Let us choose $\d>0$ such that
\be
\begin{cases}
1-\mathcal{S}\d>0,\\
1-\left(\b +\frac{\mathcal{M\d}}{1-\mathcal{S}\d}\right)>0.
\end{cases}
\ee
Therefore, since as $\e \rightarrow 0$
$$
\Big[G'_\e(u_i)-\left(\b+\frac{\mathcal{M}\d}{(1-\mathcal{S}
\d)}\right) \frac{G_\e(u_i)}{u_i} \Big]\to \Big(
1-\b-\frac{\mathcal{M}\d}{(1-\mathcal{S} \d)}\Big)>0,\,\, \quad\mbox{in
} \{u_i \neq 0\},
$$
by Fatou's Lemma we get
\begin{equation}\label{ffftatatatta}
\int_{\O\setminus \{u_i=0\}} \frac{|\n
u|^{p-2}}{|y-x|^\g}\frac{|\n u_i|^2}{|u_i|^\b}\psi dx\leq\mathcal{C}.
\end{equation}
To prove \eqref{eq:Tesi 2} we choose $E \subset \subset \O$ such
that $$ Z_u\cap (\O \setminus E)= \emptyset.$$ Since $u$ is $C^2$ in
$\O\setminus E$, then we may reduce to prove that that
$$
\int_{E\setminus \{u_i=0\}} \frac{|\n
u|^{p-2}}{|y-x|^\g}\frac{|\n u_i|^2}{|u_i|^\b} dx\leq\mathcal{C}.
$$
This, and hence the assertion, follows by considering \eqref{ffftatatatta} with a cut-off function as above with
$\psi \in C^\infty_c(\O)$ positive, such that the support of $\psi$ is compactly contained in $\O$,
$\psi \geq 0$ in $\O$ and $\psi \equiv 1$ in $E$. The proof is now complete.
\end{proof}

\subsection{Summability of $|\n u|^{-1}$}

We have the following

\begin{theorem}\label{hjfbjshdjshvb}
Let $u$ be a solution of \eqref{eq:pe} and assume, furthermore, that $f(s)>0$ for any $s>0$.
Then, there exists a positive constant $C$, independent of $y$, such that
\begin{equation}
\int_{\Omega} \frac{1}{|\nabla u|^{(p-1)r}}
\frac{1}{|x-y|^{\gamma}}\,dx\leq C
\end{equation}
where $0<r<1$ and $\gamma<n-2$ for $n\geq 3$ ($\gamma =0$ if $n=2$).

In particular the critical set $Z_u$ has zero Lebesgue measure.
\end{theorem}
\begin{proof}
Let $E$ be a  set with $E \subset\subset \Omega$ and
$(\Omega\setminus E)\cap Z_u=\emptyset$.
Recall that $Z_u=\{\nabla
u=0\}$ and $Z_u\cap \partial\Omega=\emptyset$, in view of Hopf
boundary lemma (see~\cite{PSB}). It is easy to see that, to prove the
result, we may reduce to show that
\begin{equation}
\int_{E} \frac{1}{|\nabla u|^{(p-1)r}}
\frac{1}{|x-y|^{\gamma}}\,dx\leq C
\end{equation}
To achieve this, let us consider the function
\begin{equation}
\Psi(x)=\Psi_{\varepsilon,y}(x)=\frac{1}{(|\nabla
u|+\varepsilon)^{(p-1)r}} \frac{1}{|x-y|^{\gamma}}\varphi,
\end{equation}
where $0<r<1$ and $\gamma<n-2$ for $n\geq 3$ ($\gamma =0$ if $n=2$). We also
assume that $\varphi$ is a positive $C^\infty_c(\O)$ cut-off function
such that $\varphi\equiv 1$ in $E$. Using $\Psi$ as test function in~\eqref{eq:pe},
since $f(u) \geq \sigma$ for some $\sigma>0$ in the support of $\Psi$, we get
\begin{equation*}
\begin{split}
& \s\int_\Omega \,\Psi\,dx\leq \int_{\Omega}f(u)\Psi\,dx=\int_\Omega a(u)|\nabla u|^{p-2}(\nabla u,\nabla\Psi)+ \frac 1p a'(u)|\nabla u|^p\Psi\,dx\\
&\leq  \int_\Omega a(u)|\nabla u|^{p-2}|(\nabla u,\nabla |\nabla u|)|\frac{1}{(|\nabla u|+\varepsilon)^{(p-1)r+1}} \frac{1}{|x-y|^{\gamma}}\varphi\,dx\\
&+ \int_\Omega a(u)|\nabla u|^{p-2}|(\nabla u,\nabla \frac{1}{|x-y|^{\gamma}})|\frac{1}{(|\nabla u|+\varepsilon)^{(p-1)r}} \varphi\,dx\\
&+\int_\Omega a(u)|\nabla u|^{p-2}|(\nabla u,\nabla \varphi)|\frac{1}{(|\nabla u|+\varepsilon)^{(p-1)r}} \frac{1}{|x-y|^{\gamma}}\,dx\\
&+ \int_\Omega \frac{a'(u)}{p}|\nabla u|^p\frac{1}{(|\nabla
u|+\varepsilon)^{(p-1)r}}
\frac{1}{|x-y|^{\gamma}} \varphi\,dx.
\end{split}
\end{equation*}
Consequently, we have
\begin{equation*}
\begin{split}
&\int_\Omega \Psi\,dx\leq C \bigg [ \int_\Omega |\nabla u|^{p-1} |D^2 u|\frac{1}{(|\nabla u|+\varepsilon)^{(p-1)r+1}} \frac{1}{|x-y|^{\gamma}}\varphi\,dx\\
&+ \int_\Omega \frac{|\nabla u|^{p-1}}{(|\nabla
u|+\varepsilon)^{(p-1)r}}\frac{1}{|x-y|^{\gamma+1}}
\varphi\,dx\\
&+\int_\Omega \frac{|\nabla u|^{p-1}}{(|\nabla u|+\varepsilon)^{(p-1)r}} \frac{1}{|x-y|^{\gamma}}\,dx\\
&+ \int_\Omega \frac{|\nabla u|^{p}}{(|\nabla
u|+\varepsilon)^{(p-1)r}} \frac{1}{|x-y|^{\gamma}}\,dx \bigg ].
\end{split}
\end{equation*}
Then, denoting by $C_i$, suitable positive
constants independent of $y$ and by $C_\delta$ a positive constant
depending on $\d$, we obtain
\begin{equation}
	\label{eq:esitmatesGrad1}
\begin{split}
&\int_\Omega \,\Psi\,dx\leq C_1 \int_\Omega |\nabla u|^{p-1} |D^2 u|\cdot\frac{1}{(|\nabla u|+\varepsilon)^{(p-1)r+1}}\cdot \frac{1}{|x-y|^{\gamma}}\cdot\varphi\,dx\\
&+ C_2 \int_\Omega \frac{1}{|x-y|^{\gamma+1}} \,dx + C_3 \int_\Omega  \frac{1}{|x-y|^{\gamma}}\,dx\\
&\leq  C_1 \int_\Omega |\nabla u|^{p-1} |D^2 u|\cdot\frac{1}{(|\nabla u|+\varepsilon)^{(p-1)r+1}}\cdot \frac{1}{|x-y|^{\gamma}}\cdot\varphi\,dx+C_4\\
&\leq \delta C_5\int_\Omega \frac{1}{(|\nabla u|+\varepsilon)^{(p-1)r}}\cdot \frac{1}{|x-y|^{\gamma}}\cdot\varphi\,dx\\
&+ C_\delta\int_\Omega |\nabla u|^{(p-2)-(p(r-1)+2-r)}
|D^2 u|^2\cdot \frac{1}{|x-y|^{\gamma}}\cdot\varphi\,dx
+C_6\leq\\
&\leq C_5\delta\int_\Omega \,\Psi\,dx +C_\delta.
\end{split}
\end{equation}
Here we have we used that $u \in C^{1, \a}(\O)$, $\g < n-2$ and we
have exploited the regularity result of Proposition~\ref{pr:Summability}.
Then, by \eqref{eq:esitmatesGrad1}, fixing $\delta$ sufficiently
small, such that $1-C_5\delta>0$, one concludes
\begin{equation}
\int_\Omega \,\frac{1}{(|\nabla
u|+\varepsilon)^{(p-1)r}} \frac{1}{|x-y|^{\gamma}}\varphi\,dx\leq K,
\end{equation}
for some positive constant $K$ independent of $y$. Taking the limit for
$\varepsilon$ going to zero, the assertion immediately follows by Fatou's Lemma.
\end{proof}

Proposition~\ref{hjfbjshdjshvb} provides in fact the right summability
of the weight $\rho(x)=|\nabla u(x)|^{p-2}$ in order to obtain a weighted Poincar\'e inequality.
We refer the readers to \cite[Section 3]{DS1} for further details. For the sake of selfcontainedness,
we recall here the statement

\begin{theorem}\label{T:CRU}
If $u \in C^{1, \a}(\overline{\Omega})$ is a solution of
\eqref{eq:pe} with $f(s) >0$ for $s>0$, $p\geq 2$, then
\begin{equation}\label{East:3}
\|v\|_{L^q(\Omega)} \leq C_p(|\Omega|)\|\nabla
v\|_{L^q(\Omega,\rho)}, \qquad \text{for every $v \in
H^{1,q}_{0,\rho}(\Omega)$},
\end{equation}
where $\rho \equiv |\nabla u|^{p-2}$, $C_P(|\Omega|) \rightarrow 0$
if $|\Omega| \rightarrow 0$. In particular \eqref{East:3} holds
for every function $ v \in H^{1,2}_{0,\rho}(\Omega) $.
Moreover if $p \geq 2$, $q \geq 2$ and $v \in
W^{1,q}_0(\Omega)$, the same conclusion holds. In fact, being $u \in
C^{1, \a}(\overline{\Omega})$, and  $p \geq 2 $, $\rho =
|Du|^{p-2}$ is bounded, so that $W^{1,q}_0(\Omega) \hookrightarrow
H^{1,q}_{0,\rho}(\Omega)$.
\end{theorem}

Recall that, if $\rho \in L^{1}(\Omega)$, $1\leq q<\infty$, the
space $H^{1,q}_\rho(\Omega)$ is defined
as the completion of $C^1(\overline{\Omega})$ (or $C^{\infty }(\overline{\Omega})$) under the norm
  \begin{equation}\label{hthInorI}
\| v\|_{H^{1,q}_\rho}= \| v\|_{L^q (\Omega)}+\| \nabla v\|_{L^q
(\Omega, \rho)}
\end{equation}
where
$$
\| \nabla v\|^q_{L^q (\Omega, \rho)}=\int_{\Omega}|\nabla v|^q \rho\,dx.
$$
 We also recall that $H^{1,q}_{0, \rho}$ may be equivalently defined as the space of
functions having  distributional derivatives represented by a
function for which the norm defined in~\eqref{hthInorI} is bounded.
These two definitions are equivalent if the domain has piecewise
regular boundary (as we are indeed assuming).

\subsection{Comparison principles}

We now have the following

\begin{proposition}
	\label{pro:confrbis}
Let $\tilde{\O}$ be a bounded smooth domain such that
$\tilde{\O}\subseteq \O.$ Assume that $u,v$  are solutions to
the problem \eqref{eq:pe} and assume that $u\leq v$ on $\partial
{\tilde \O}$. Then there exists a positive constant $\theta$, depending both on $u$ and $f$,
such that, assuming
$$
\mathcal{L}(\tilde{\O}) \leq \theta
$$
then it holds
\begin{equation*}
u\leq v\quad \text{ in  ${\tilde\O}$}.
\end{equation*}
\end{proposition}

\begin{proof}
We start proving the result when $p>2$.
Let us recall the weak formulations
\begin{align}
	\label{eqA}
& \int_{\Omega}a(u)|\nabla u|^{p-2}(\nabla u\, ,\, \nabla
\varphi)\,+\, \frac{a'(u)}{p}|\nabla u|^p\varphi\,dx=\int_\Omega
f(u)\varphi \,dx,   \\
	\label{eqB}
& \int_{\Omega}a(v)|\nabla v|^{p-2}(\nabla v\, ,\, \nabla
\varphi)\,+\, \frac{a'(v)}{p}|\nabla v|^p\varphi\,dx=\int_\Omega
f(v)\varphi \,dx.
\end{align}
Then we assume by contradiction that the assertion is false, and consider
$$
(u-v)^+=\max\{u-v,0\},
$$
that, consequently, is not identically equal to zero. Let us also set
$\Omega ^+\equiv {\rm supp} (u-v)^+\cap \wt{\O}$.
Since by assumption $u\leq v$ on $\partial \wt{\O}$,
it follows that $(u-v)^+\in W^{1,p}_0(\wt{\O})$. We can therefore
choose it as admissible test function in~\eqref{eqA} and~\eqref{eqB}. Whence,
subtracting the two, we get
\begin{equation}
\begin{split}
&\int_{\Omega^+}a(u)|\nabla u|^{p-2}(\nabla u\, ,\, \nabla (u-v))-a(v)|\nabla v|^{p-2}(\nabla v\, ,\, \nabla (u-v))\,+\\
&+\int_{\Omega^+} \frac{a'(u)}{p}|\nabla u|^p(u-v)\,dx\, -\frac{a'(v)}{p}|\nabla v|^p(u-v)\,dx=\\
&=\int_{\Omega^+} (f(u)-f(v))(u-v)\,dx.
\end{split}
\end{equation}
We can rewrite as follows
\begin{equation}
\begin{split}
&\int_{\Omega^+}a(u)((|\nabla u|^{p-2}\nabla u-|\nabla v|^{p-2}\nabla v)\, ,\, \nabla (u-v)))\,dx\\
&+\int_{\Omega^+} (a(u)-a(v))|\nabla v|^{p-2}(\nabla v,\nabla (u-v))dx\\
&+\int_{\Omega^+} \frac{1}{p}(a'(u)-a'(v))|\nabla u|^p(u-v)\,dx\, \\
&+\int_{\Omega^+} \frac{a'(v)}{p}(|\nabla u|^p-|\nabla v|^p)(u-v)\,dx \\
&=\int_{\Omega^+} (f(u)-f(v))(u-v)\,dx.
\end{split}
\end{equation}
First of all, since $a(u)\geq \eta>0$, and using the fact that
$$
\left (|\xi|^{p-2}\xi-|\xi'|^{p-2}\xi',\xi-\xi'\right)\geq c (|\xi|+|\xi'|)^{p-2}|\xi-\xi'|^2
$$
for all $\xi,\xi' \in \R^n$, it follows that
\begin{equation}
\begin{split}
c\eta \int_{\Omega^+}(|\nabla u|+|\nabla v|)^{p-2}|\nabla
(u-v)|^2\,dx\leq
&\int_{\Omega^+}a(u)(|\nabla u|^{p-2}\nabla u-|\nabla v|^{p-2}\nabla v , \nabla (u-v))\,dx\\
\end{split}
\end{equation}
so that
\begin{equation}\label{ciccipizzi}
\begin{split}
 &\int_{\Omega^+}(|\nabla u|+|\nabla v|)^{p-2}|\nabla (u-v)|^2\,dx\leq C\int_{\Omega^+} |a(u)-a(v)||\nabla v|^{p-1}|\nabla (u-v)|dx+\\
&+C\int_{\Omega^+} |a'(u)-a'(v)||\nabla u|^p|u-v|\,dx\, +C\int_{\Omega^+} |a'(v)||\nabla u|^p-|\nabla v|^p||u-v|\,dx+\\
&+\int_{\Omega^+} |\frac{f(u)-f(v)}{u-v}||u-v|^2\,dx
\end{split}
\end{equation}
Let us now evaluate the terms on right of the above inequality.
By the smoothness of $a$, the $C^{1,\alpha}$ regularity of
$u$, and exploiting Young inequality we get
\begin{equation}
\begin{split}
 &\int_{\Omega^+} |a(u)-a(v)||\nabla v|^{p-1}|\nabla (u-v)|dx\leq C\int_{\Omega^+} |u-v||\nabla v|^{\frac{p-2}{2}}|\nabla (u-v)|dx\leq \\
 &\leq C_\delta\int_{\Omega^+} (u-v)^2\,dx+\delta\int_{\Omega^+} (|\nabla u|+|\nabla v|)^{p-2}|\nabla (u-v)|^2\,dx\leq\\
 &\leq (C_\delta C_p(|\Omega^+|)+\delta)\int_{\Omega^+} (|\nabla u|+|\nabla v|)^{p-2}|\nabla (u-v)|^2\,dx.
\end{split}
\end{equation}
Here $C_\delta$ is a constant depending on $\delta$, and
$C_p(|\Omega^+|)$ is the Poincar\'e constant given by Theorem~\ref{T:CRU}.
Note in particular that, since $p>2$, we have $|\nabla u|^{p-2}\leq(|\nabla u|+|\nabla v|)^{p-2}$.
It is of course very important the fact that the constant $C_p(|\Omega^+|)$
goes to zero, provided that the Lebesgue measure of $\Omega^+$ goes to $0$.
Also we note that, by the $C^{1,\alpha}$ regularity of $u$, and exploiting
the fact that $a'$ is Lipschitz continuous, we get
\begin{align*}
\int_{\Omega^+} |a'(u)-a'(v)||\nabla u|^p|u-v|\,dx &\leq
C\int_{\Omega^+}(u-v)^2\,dx   \\
&\leq
C\,C_P(|\Omega^+|)\int_{\Omega^+}(|\nabla u|+|\nabla
v|)^{p-2}|\nabla (u-v)|^2\,dx.
\end{align*}
Also, by convexity, we have
\begin{equation}
\begin{split}
&\int_{\Omega^+} |a'(v)||\nabla u|^p-|\nabla v|^p||u-v|\,dx \\
&\leq C\int_{\Omega^+} (|\nabla u|+|\nabla v|)^{\frac{p-2}{2}}|\nabla (u-v)||u-v|\,dx  \\
&\leq \delta\int_{\Omega^+} (|\nabla u|+|\nabla v|)^{p-2}|\nabla (u-v)|^2\,dx+C_\delta\int_{\Omega^+} |u-v|^2\,dx \\
&\leq \delta\int_{\Omega^+} (|\nabla u|+|\nabla v|)^{p-2}|\nabla (u-v)|^2\,dx \\
&+C_\delta C_P(|\Omega ^+|)\int_{\Omega^+} (|\nabla u|+|\nabla v|)^{p-2}|\nabla (u-v)|^2\,dx  \\
&\leq (\delta+C_\delta C_P(|\Omega ^+|))\int_{\Omega^+} (|\nabla u|+|\nabla v|)^{p-2}|\nabla (u-v)|^2\,dx
\end{split}
\end{equation}
Finally, by the Lipschitz continuity of $f$, it follows
\begin{align*}
\int_{\Omega^+} |\frac{f(u)-f(v)}{u-v}||u-v|\,dx&\leq
C\int_{\Omega^+} |u-v|^2\,dx \\
&\leq C\,C_P(|\Omega ^+|)\int_{\Omega^+} (|\nabla u|+|\nabla
v|)^{p-2}|\nabla (u-v)|^2\,dx
\end{align*}
Concluding, exploiting the above estimates, we get
\begin{equation*}
\int_{\Omega^+} (|\nabla u|+|\nabla v|)^{p-2}|\nabla
(u-v)|^2\,dx\leq (\delta +C_\delta\, C_P(|\Omega ^+|))\int_{\Omega^+} (|\nabla
u|+|\nabla v|)^{p-2}|\nabla (u-v)|^2\,dx
\end{equation*}
which gives a contradiction for
$(\delta +C_\delta\, C_P(|\Omega ^+|))<1$. Therefore, if we consider $\delta $ small fixed, say $\delta=\frac{1}{4}$, it then follows that also $C_\delta$ is fixed.
Now, since $\mathcal{L}(\wt{\O}) \leq \theta$ by assumption, it follows that if $\theta $ is sufficiently small,
then we may assume that $C_P(|\Omega ^+|)$ is also small, and that  $C_\delta\, C_P(|\Omega ^+|))<\frac{1}{4}$. Consequently, it follows $(\delta +C_\delta\, C_P(|\Omega ^+|))<\frac{1}{2}<1$, that leads to the above contradiction, and shows that actually $(u-v)^+=0$ and the thesis.
The proof in the case $1<p\leq 2$ in completely analogous, but is based on the classical Poincar\'e inequality.
We give some details for the reader's convenience.
Exactly as above we get \eqref{ciccipizzi}. This , for $1<p\leq 2$, considering the fact that the term
$(|\nabla u|+|\nabla v|)^{p-2}$ is bounded below by the fact that $p-2\leq 0$ and $|\nabla u|\, ,\, |\nabla v|\in L^\infty(\overline{\Omega})$, gives
\begin{equation}\label{ciccipizzi2}
\begin{split}
 &\int_{\Omega^+}|\nabla (u-v)|^2\,dx\leq C\int_{\Omega^+} |a(u)-a(v)||\nabla v|^{p-1}|\nabla (u-v)|dx+\\
&+C\int_{\Omega^+} |a'(u)-a'(v)||\nabla u|^p|u-v|\,dx\, +C\int_{\Omega^+} |a'(v)|||\nabla u|^p-|\nabla v|^p||u-v|\,dx+\\
&+\int_{\Omega^+} |\frac{(f(u)-f(v))}{(u-v)}|\cdot||u-v|\,dx\leq \\
&C\int_{\Omega^+} |u-v||\nabla (u-v)|dx+C\int_{\Omega^+} |u-v|^2\,dx\, \leq \\
&\delta\int_{\Omega^+} |\nabla (u-v)|^2dx+ C_\delta\int_{\Omega^+} |u-v|^2\,dx\, \leq\\
&\delta\int_{\Omega^+} |\nabla (u-v)|^2dx+C_\delta C_P(|\Omega^+|)\int_{\Omega^+} |\nabla (u-v)|^2\,dx\, \leq\\
&(\delta+C_\delta
C_P(|\Omega^+|))\int_{\Omega^+} |\nabla (u-v)|^2\,dx
\end{split}
\end{equation}
For  $\theta $ sufficiently small arguing as above we can assume $(\delta+C_\delta
C_P(|\Omega^+|))<1$ which
gives $(u-v)^+=0$ and the thesis.
\end{proof}

\subsection{The moving plane method}
Let us consider a direction, say $x_1$, for example. As customary we set
$$
T_\lambda=\big\{x \in \R^n: x_1= \lambda \big \}.
$$
Given $x\in \R^n$, we define
$$
x_\lambda= (2\lambda-x_1, x_2,\ldots,x_n), \quad
u_\lambda(x)=u(x_\lambda),
$$
$$
\O_\lambda=\big \{x\in \O: x_1<\lambda \big\},
$$
Set
$$
\tilde{a}:=\inf_{x\in \O} x_1.
$$
Let $\Lambda $ be the set of those  $\lambda >\tilde{a}$ such that for each $\mu < \lambda$
none of the conditions (i) and (ii) occurs, where
 \begin{itemize}
 \item[(i)] The reflection of $(\Omega_{\lambda})$ w.r.t. $T_\lambda$ becomes internally tangent to $\partial \Omega$ .
 \item[(ii)]$T_{\lambda}$ is orthogonal to $\partial \Omega$.
 \end{itemize}

We have the following

\begin{proposition}\label{prosimmetri2}
Let $u\in C^{1,\a}(\ov{\O})$ be a solution to the problem \eqref{eq:pe}.
Then, for any $\tilde{a}\leq\lambda\leq \Lambda$, we have
\begin{equation}\label{zxclu:1.8222}
u(x)\leq u_{\lambda}(x),\qquad \forall x \in \Omega _{\lambda}.
\end{equation}
Moreover, for any $\lambda$ with $\tilde{a}<\lambda <\Lambda$ we have
\begin{equation}\label{zxcNEVFURu}
u(x) < u_{\lambda}(x),\qquad \forall x \in \Omega_{\lambda}\setminus
Z_{u,\lambda}
\end{equation}
where $Z_{u,\lambda}\equiv\{x\in \Omega _{\lambda}\, :\, \nabla
u(x)=\nabla u_{\lambda}(x)=0\}$. Finally
\begin{equation}\label{zxclu:1.9222u}
\frac{\partial u}{\partial x_1}(x) \geq 0, \qquad \forall x \in
\Omega _{\Lambda}.
\end{equation}
\end{proposition}

\begin{proof}
For $\tilde{a}<\lambda<\Lambda$ and $\lambda$ sufficiently close to $\tilde{a}$, we
assume that
 $\mathcal{L}(\O_\lambda)$ is as small as we like.
 We assume in particular that we can exploit the weak maximum principle in small domains (see Proposition \ref{pro:confrbis}) in $\O_\lambda$.
 Consequently, since we know that
\begin{equation}
u-u_{\lambda} \leq 0,\qquad\text{on $\partial\O_\lambda$}
\end{equation}
 by construction, by Proposition \ref{pro:confrbis} it follows
\begin{eqnarray*}
u-u_{\lambda} \leq 0 \quad\mbox{in}\quad \O_\lambda.
\end{eqnarray*}
We define
\begin{equation}
\Lambda_0=\{\lambda > \tilde{a}: u\leq u_{t},\,\,\text{for all
$t\in(\tilde{a},\lambda]$}\}
 \end{equation}
and
 \begin{equation}
\lambda_0=\sup\,\Lambda_0.
\end{equation}
Note that by continuity, we have $u\leq u_{\lambda_0}$.
We have to show that actually $\lambda_0 = \Lambda$. Assume that by
contradiction $\lambda_0 <\Lambda$ and argue as follows.
Let  $A$ be an open set such that $Z_u\cap\Omega_{\lambda_0} \subset
A \subset \Omega_{\lambda_0}$. Note that since $|Z_u|=0$
(see Theorem \ref{hjfbjshdjshvb}), we can choice $A$ as small as we like.
Note now that by a strong comparison principle \cite{PSB}
we get
$$
u<u_{\lambda_0} \qquad\text{or}\qquad u\equiv u_{\lambda_0}
$$
in any connected component of $\Omega_{\lambda_{0}}\setminus Z_u$.\\
It follows now that
\begin{center}
the case $u\equiv u_{\lambda_0}$ in some  connected component
$\mathcal{C}$ of $\Omega_{\lambda_{0}}\setminus Z_u$ is not
possible.
\end{center}
The proof of this is completely analogous to the one given in \cite{DP} once we have Proposition \ref{pro:confrbis}.
Consider now a compact set $K$ in $\Omega_{\lambda_{0}}$ such that
$|\Omega_{\lambda_{0}}\setminus K|$ is sufficiently small so that
Proposition \ref{pro:confrbis} works. By what we proved before, $u_{\lambda_{0}}-u$ is positive in $K
\setminus A$ which is compact, therefore  by continuity we find
$\epsilon> 0$ such that, $\lambda_0+\epsilon < \Lambda$ and for $
\lambda<\lambda_0+\epsilon$ we have that $|\Omega_{\lambda}\setminus
(K\setminus A)|$ is still sufficiently small as before and
$u_{\lambda}-u >0$ in $K \setminus A$. In particular $u_{\lambda}-u
>0$ on $\partial(K \setminus A)$. Consequently $u \leq
u_\lambda$  on $\partial(\Omega_\lambda\setminus(K\setminus A))$. By
Proposition \ref{pro:confrbis} it follows $u \leq u_\lambda$ in
$\Omega_\lambda\setminus(K\setminus A)$ and consequently in
$\Omega_\lambda$, which contradicts the assumption $\lambda_0 <
\Lambda$.
Therefore $\lambda_0 \equiv \Lambda$ and the thesis is proved.
The proof of \eqref{zxcNEVFURu} follows by the strong comparison theorem exploited as above.
Finally \eqref{zxclu:1.9222u}  follow by the monotonicity of the
solution that is implicitly in the above arguments.
\end{proof}

\medskip

\section{Properties of the parabolic flow}
\label{parabolic}

Let $\Omega$ be a smooth bounded domain in $\R^n$, and let $a:\R\to\R$ be a $C^1$ function
such that there exists positive constants $C,\nu$ and $\rho$ such that
\begin{align}
	\label{boundd}
& \eta\leq a(s)\leq C,\,\,\, |a'(s)|\leq C\quad\text{for all $s\in\R$}, \\
& a'(s)s\geq 0,\quad\text{for all $s\in\R$ with $|s|\geq \rho$}.
\label{sign}
\end{align}
As stated in the introduction,
along any given global solution $u:\R^+\times\Omega\to\R$ of problem~\eqref{problema}, and setting
$$
F(s)=\int_0^s f(\tau)d\tau,\quad s\in\R,
$$
we also consider the energy functional ${\mathcal E}$ defined by
$$
{\mathcal E}(u(t))=\frac{1}{p}\int_{\Omega}a(u(t))|\nabla u(t)|^pdx-
\int_{\Omega}F(u(t))dx,
$$
and the related energy inequality~\eqref{energyinequality}.
In particular, the energy functional ${\mathcal E}$ is
non-increasing along solutions. Moreover, by the regularity we assumed
on the global solutions, we have
\begin{equation}
	\label{boundedtj}
\sup_{t>0}\|u(t)\|_{W^{1,p}_0(\O)}<\infty,
\end{equation}
and
\begin{equation}
	\label{summabilitygg}
	\int_0^\infty\int_{\O}|u_t(\tau)|^2dxd\tau<\infty.
\end{equation}
\vskip3pt
\noindent
Next we state a quite useful result.

\begin{lemma}
	\label{elle2conv}
	For all fixed $\mu_0>0$, it holds
	$$
	\lim_{t\to\infty}\sup_{\mu\in [0,\mu_0]}\|u(t)-u(t+\mu)\|_{L^q(\Omega)}=0,
	\quad\text{for all $q\in [1,p^*)$}.
	$$
	If in addition the trajectory $\{u(t):t>1\}$ is relatively compact in $W^{1,p}_0(\Omega)$, we have
	$$
	\lim_{t\to\infty}\sup_{\mu\in [0,\mu_0]}\|u(t)-u(t+\mu)\|_{W^{1,p}_0(\Omega)}=0,
	$$
	for all fixed $\mu_0>0$.
\end{lemma}
\begin{proof}
	Let us first prove that, for all $\mu_0>0$, it holds
	\begin{equation}
		\label{elle1lim-O}
	\lim_{t\to\infty}\sup_{\mu\in [0,\mu_0]}\|u(t)-u(t+\mu)\|_{L^1(\Omega)}=0.
\end{equation}
	Given $\mu>0$, for all $t>0$ and $\mu\in [0,\mu_0]$, from the energy inequality~\eqref{energyinequality},
	we have
	\begin{align*}
	\int_{\Omega}|u(t)-u(t+\mu)|dx & =\int_{\Omega}\Big|\int_t^{t+\mu}u_t(\tau)d\tau\Big|dx
	\leq  \int_t^{t+\mu}\int_{\Omega} |u_t(\tau)|d\tau dx \\
	& \leq \sqrt{\mu{\mathcal L}^n(\Omega)} \Big(\int_t^{t+\mu}\int_{\Omega}|u_t(\tau)|^2d\tau dx\Big)^{1/2} \\
	\noalign{\vskip3pt}
&	\leq \sqrt{\mu{\mathcal L}^n(\Omega)} ({\mathcal E}(u(t))-{\mathcal E}(u(t+\mu)))^{1/2} \\
	\noalign{\vskip5pt}
&	\leq \sqrt{\mu_0{\mathcal L}^n(\Omega)} ({\mathcal E}(u(t))-{\mathcal E}(u(t+\mu_0)))^{1/2}.
\end{align*}
	Then, since ${\mathcal E}$ is non-increasing and bounded below, the assertion follows by letting
	$t\to\infty$ in the previous inequality. Let now $q\in [1,p^*)$ and assume now by contradiction that along
	a diverging sequence of times $(t_j)$, we get
	$$
	\sup_{\mu\in [0,\mu_0]}\|u(t_j)-u(t_j+\mu)\|_{L^q(\Omega)}\geq\sigma>0,
	$$
	for some positive constant $\sigma$ and all $j$ large. In particular, there is a sequence $\mu_j\subset
	[0,\mu_0]$ such that $\|u(t_j)-u(t_j+\mu_j)\|_{L^q(\Omega)}\geq\sigma>0$ for all $j$ large.
	In light of~\eqref{boundedtj}, by Rellich compactness Theorem, up to a subsequence, it follows that $u(t_j)\to\xi_1$
	in $L^q(\Omega)$ as $j\to\infty$ and $u(t_j+\mu_j)\to\xi_2$
	in $L^q(\Omega)$ as $j\to\infty$, yielding $\|\xi_2-\xi_1\|_{L^q(\Omega)}\geq\sigma>0$. In particular
	$\xi_1\neq\xi_2$. On the other hand, from~\eqref{elle1lim-O} we immediately get
	$\|\xi_2-\xi_1\|_{L^1}=0$, leading to a contradiction. The second part of the statement has an analogous
	proof assuming by contradiction that there exists $\sigma>0$ and a diverging sequence of times $(t_j)$
	such that
	$$
	\sup_{\mu\in [0,\mu_0]}\|u(t_j)-u(t_j+\mu)\|_{W^{1,p}_0(\Omega)}\geq\sigma>0,
	$$
	and then exploiting the relative compactness
	of $\{u(t):t>1\}$ in $W^{1,p}_0(\Omega)$.
\end{proof}

On $W^{1,p}_0(\Omega)$ the functional ${\mathcal E}$ is defined by setting
\begin{equation}
\label{generalmod}
{\mathcal E}(u)=\frac{1}{p}\int_\Omega a(u)|\nabla u|^p-\int_\Omega F(u).
\end{equation}
and it is merely continuous, although its directional
derivatives exist along smooth directions and
\begin{equation*}
{\mathcal E}'(u)(\varphi)=\int_\Omega a(u)|\nabla u|^{p-2}\nabla u\cdot\nabla\varphi
+\frac{1}{p}\int_\Omega a'(u)|\nabla u|^p\varphi
-\int_\Omega f(u)\varphi.
\end{equation*}

We now recall an important compactness result (see e.g.~\cite{CD,Sq1}).

\begin{lemma}
	\label{strongbound}
	Let conditions~\eqref{boundd} and~\eqref{sign} hold.
Assume that $(u_h)\subset W^{1,p}_0(\Omega)$ is a bounded sequence and
\begin{equation*}
\langle w_h,\varphi\rangle=
\int_{\Omega}a(u_h)|\nabla u_h|^{p-2}\nabla u_h\cdot\nabla\varphi+
\frac{1}{p}\int_{\Omega}a'(u_h)|\nabla u_h|^p\varphi
\end{equation*}
for every $\varphi\in C^\infty_c(\Omega)$,
where $(w_h)$ is strongly convergent in $W^{-1,p'}(\Omega)$.
Then $(u_h)$ admits a strongly convergent subsequence in $W^{1,p}_0(\Omega)$.
\end{lemma}

\begin{lemma}
	\label{parab-elliptic-th}
	Let conditions~\eqref{boundd} and~\eqref{sign} hold.
	Assume that there exist $C_1,C_2>0$ such that
	\begin{equation}
		\label{growthf}
	|f(s)|\leq C_1+C_2|s|^r,\qquad\text{for all $s\in\R$},
\end{equation}
	for some $r\in [1,p^*-1)$.
Let $u:[0,\infty)\times\Omega\to\R$ be a global solution to problem~\eqref{problema}, with $p>\frac{2n}{n+2}$.
Then, for every diverging sequence $(\tau_j)$ there exists a diverging
sequence $(t_j)$ with $t_j\in [\tau_j,\tau_j+1]$ such that
\begin{equation}
	\label{convinfo}
u(t_j)\to z\quad\text{in $W^{1,p}_0(\O)$ as $j\to\infty$},
\end{equation}
where either $z=0$ or $z$ is a solution to problem~\eqref{eq:pe}. In addition, it holds
	$$
	\lim_{t\to\infty}\sup_{\mu\in [0,\mu_0]}\|u(t_j+\mu)-z\|_{L^q(\Omega)}=0,
	\quad\text{for all $q\in [1,p^*)$},
	$$
for all fixed $\mu_0>0$.
\end{lemma}
\begin{proof}
	By the definition of solution, for all
	$\varphi\in C^\infty_c(\Omega)$ and for a.e. $t>0$, we have
	\begin{align}
		\label{defweakae}
	\int_{\Omega} u_t(t)\varphi dx &+\int_{\Omega}a(u(t))|\nabla u(t)|^{p-2}\nabla u(t)\cdot\nabla \varphi dx  \\
&	+\int_{\Omega}\frac{a'(u(t))}{p}|\nabla u(t)|^{p}\varphi dx =\int_{\Omega}f(u(t))\varphi dx.  \notag
\end{align}
	By means of the summability given by~\eqref{summabilitygg} it follows that, for every diverging sequence
	$(\tau_j)\subset\R^+$, there exists a diverging sequence $(t_j)$ with $t_j\in [\tau_j,\tau_j+1]$, $j\geq 1$, such that
	\begin{equation}
		\label{vanishingSQ}
	\Lambda_j=\int_\Omega |u_t(t_j)|^2dx\to 0,\quad\text{as $j\to\infty$}.
\end{equation}
	Let us now define the sequence $(w_j)$ in $W^{-1,p'}(\Omega)$ by
	$$
	\langle w_j,\varphi\rangle=\langle w_j^1,\varphi\rangle+\langle w_j^2,\varphi\rangle,
	\qquad\text{for all $\varphi\in W^{1,p}_0(\O)$},
	$$
	where we have set
	$$
	\langle w_j^1,\varphi\rangle=\int_{\O}f(u(t_j))\varphi,\quad
	\langle w_j^2,\varphi\rangle=-\int_{\O} u_t(t_j)\varphi\,dx,
	\qquad\text{for all $\varphi\in W^{1,p}_0(\O)$}.
	$$
	We recall that, under the growth condition~\eqref{growthf}, the map
	$$
	W^{1,p}_0(\O)\ni u\mapsto f(u)\in W^{-1,p'}(\O)
	$$
	is completely continuous, and hence, up to a further subsequence, we have
	$$
	w_j^1\to \mu,\qquad\text{in $W^{-1,p'}(\O)$ as $j\to\infty$},
	$$
	for some $\mu\in W^{-1,p'}(\O)$. Turning to the sequence $(w_j^2)$, notice that
	in view of~\eqref{vanishingSQ}, exploiting the fact that $p^*>2$ since of the assumption $p>\frac{2n}{n+2}$,
by H\"older inequality we get
	$$
	\|w_j^2\|_{W^{-1,p'}(\Omega)}=\sup\big\{|\langle w_j,\varphi  \rangle|: \varphi\in W^{1,p}_0(\O),\,\,
	\|\varphi\|_{W^{1,p}_0(\O)}\leq 1\big\}\leq C\Lambda_j,
	$$
	for some positive constant $C$. Then $w_j^2\to 0$ in $W^{-1,p'}(\Omega)$ as $j\to\infty$
	and, in conclusion, $w_j\to \mu$ in $W^{-1,p'}(\Omega)$ as $j\to\infty$. Furthermore,
	by means of~\eqref{defweakae}, we conclude that
	\begin{equation}
		\label{variationalseq}
\langle w_j,\varphi\rangle=
\int_{\Omega}a(u(t_j))|\nabla u(t_j)|^{p-2}\nabla u(t_j)\cdot\nabla\varphi+
\frac{1}{p}\int_{\Omega}a'(u(t_j))|\nabla u(t_j)|^p\varphi,
\end{equation}
for all $\varphi\in C^\infty_c(\Omega)$.
We have thus proved that $(u(t_j))\subset W^{1,p}_0(\O)$ is in the framework
of the compactness Lemma~\ref{strongbound}.
In turn, by Lemma~\ref{strongbound}, up to a subsequence
$(u(t_j))$ is strongly convergent to some $z$ in $W^{1,p}_0(\O)$, as $j\to\infty$. In particular,
$u(t_j,x)\to z(x)$ and $\nabla u(t_j,x)\to \nabla z(x)$ for a.e. $x\in\O$, as $j\to\infty$.
Since
$$
|a'(u(t_j,x))|\nabla u(t_j,x)|^{p}\varphi(x)|\leq C|\nabla u(t_j,x)|^{p}, \quad\text{for all $j\geq 1$ and $x\in\O$},
$$
and
$|\nabla u(t_j,x)|^{p}\to |\nabla z(x)|^p$ in $L^1(\O)$ as $j\to\infty$, we have
$$
\lim_{j\to\infty}\int_{\Omega}a'(u(t_j))|\nabla u(t_j)|^{p}\varphi dx
=\int_{\Omega}a'(z)|\nabla z|^{p}\varphi dx
$$
by generalized Lebesgue dominated convergence theorem. Also, as
$$
a(u(t_j,x))|\nabla u(t_j,x)|^{p-2}\nabla u(t_j,x)\to a(z(x))|\nabla z(x)|^{p-2}\nabla z(x),
$$
and
$$
\text{$a(u(t_j))|\nabla u(t_j)|^{p-2}\nabla u(t_j)$\,\,\, is bounded in $L^{p'}(\Omega)$},
$$
we have
$$
\lim_{j\to\infty} \int_{\Omega}a(u(t_j))|\nabla u(t_j)|^{p-2}\nabla u(t_j)\cdot\nabla \varphi \,dx
=\int_{\Omega}a(z)|\nabla z|^{p-2}\nabla z\cdot\nabla \varphi \,dx
$$
Finally, since $f(u(t_j,x))\to f(z(x))$ a.e. in $\O$, as $j\to\infty$, we get
$$
\lim_{j\to\infty}\langle w_j,\varphi\rangle=
\lim_{j\to\infty}\int_{\Omega}f(u(t_j))\varphi \,dx=\int_{\Omega}f(z)\varphi \,dx.
$$
In particular, letting $j\to\infty$ in formula~\eqref{variationalseq},
it follows that $z$ is a (possibly zero) weak solution to problem
$$
-{\rm div}(a(z) |\n z|^{p-2}\n z)+\frac{a'(z)}{p}|\n z|^p = f(z),\quad \text {in $\O$}.
$$
The last assertion of the statement is just a combination of~\eqref{convinfo} with
Lemma~\ref{elle2conv}.
\end{proof}

\begin{lemma}
	\label{lemma-bis-comp}
	Let $u_0\in {\mathcal A}$ and let $u:[0,\infty)\times\Omega\to\R^+$ be the
	corresponding global solution to problem~\eqref{problema}. Then
	the $\omega$-limit set $\omega(u_0)$ only contains
	positive (possibly identically zero) solutions of problem~\eqref{eq:pe}.
\end{lemma}
\begin{proof}
Let $z\in \omega(u_0)$. Therefore, there exists a diverging sequence $(t_j)\subset\R^+$ such that
$u(t_j)$ converges to $z$ in $W^{1,p}_0(\Omega)$, as $j\to\infty$.
Let now $\varphi\in C^\infty_c(\Omega)$ be a given
test function with $\|\varphi\|_{C^1}\leq 1$. Multiply problem~\eqref{problema} by $\varphi$
and integrate it in space over $\Omega$ and in time over $[t_j,t_j+\sigma_j]$, where
$\sigma_j\in [\sigma,1]$ for a fixed $\sigma>0$, yielding
\begin{align}
		\label{allequat}
		\int_{t_j}^{t_j+\sigma_j}\int_{\Omega} u_t\varphi dx &+	\int_{t_j}^{t_j+\sigma_j}\int_{\Omega}a(u)|\nabla u|^{p-2}\nabla u\cdot\nabla \varphi dx \\
		&
	+\frac{1}{p}	\int_{t_j}^{t_j+\sigma_j}\int_{\Omega}a'(u)|\nabla u|^{p}\varphi dx=	\int_{t_j}^{t_j+\sigma_j}\int_{\Omega}f(u)\varphi dx,  \notag
	\end{align}
for any $j\geq 1$. Now, by virtue of Lemma~\ref{elle2conv}, it follows that
	\begin{align*}
\Big|		\int_{t_j}^{t_j+\sigma_j}\int_{\Omega} u_t\varphi dx\Big|&=
\Big|\int_{\Omega} (u(t_j+\sigma_j)-u(t_j))\varphi dx\Big|  \\
& \leq \int_{\Omega} |u(t_j+\sigma_j)-u(t_j)||\varphi| dx \\
\noalign{\vskip5pt}
& \leq C\|u(t_j+\sigma_j)-u(t_j)\|_{L^1}=o(1),\quad\text{as $j\to\infty$}.
\end{align*}
	In particular, recalling that $u\in C([0,\infty),W^{1,p}_0(\O,\R^+))$,
	by applying the mean value theorem, we find
	a new diverging sequence $(\xi_j)\subset\R^+$ with $\xi_j\in [t_j,t_j+\sigma_j]$ such that
		\begin{align}
			\label{quasisol}
	\int_{\Omega}a(u(\xi_j))|\nabla u(\xi_j)|^{p-2}\nabla u(\xi_j)\cdot\nabla \varphi dx
	&+\frac{1}{p}\int_{\Omega}a'(u(\xi_j))|\nabla u(\xi_j)|^{p}\varphi dx \\
&	=\int_{\Omega}f(u(\xi_j))\varphi dx
	+o(1),\quad\text{as $j\to\infty$}.  \notag
	\end{align}
	In general, the choice of the sequence $(\xi_j)$ may depend upon the
	particular test function $\varphi$ that was fixed. On the other hand, taking into account
	the second part of the statement of Lemma~\ref{elle2conv}, without loss of generality we may assume that $\xi_j$
	is independent of $\varphi$. In fact, denoting by $(\xi_j^0)$
	and $(\xi_j^\varphi)$ the sequences satisfying the property above and related
	to a reference test functions $\varphi_0$ and to an arbitrary test function $\varphi$
	respectively, and writing,
	\begin{equation}
u(\xi_j^0)-u(\xi_j^\varphi)=\beta_j,\qquad \text{where $\beta_j\to 0$ in $W^{1,p}_0(\Omega)$ as $j\to\infty$},
\end{equation}
	where $\beta_j$ is independent of $\varphi$, we get
	\begin{align*}
	&\Big|\int_{\Omega}a(u(\xi_j^0))|\nabla u(\xi_j^0)|^{p-2}\nabla u(\xi_j^0)\cdot\nabla \varphi dx-
	\int_{\Omega}a(u(\xi_j^\varphi))|\nabla u(\xi_j^\varphi)|^{p-2}\nabla u(\xi_j^\varphi)\cdot\nabla \varphi \,dx\Big|	 \\
&	=\Big|\int_{\Omega}\big(a(u(\xi_j^0))|\nabla u(\xi_j^0)|^{p-2}\nabla u(\xi_j^0)-a(u(\xi_j^\varphi))|\nabla u(\xi_j^\varphi)|^{p-2}\nabla u(\xi_j^\varphi)\big)\cdot\nabla \varphi \,dx\Big| \\
&\leq \int_{\Omega}\big|a(u(\xi_j^0))|\nabla u(\xi_j^0)|^{p-2}\nabla u(\xi_j^0)-a(u(\xi_j^\varphi))|\nabla u(\xi_j^\varphi)|^{p-2}\nabla u(\xi_j^\varphi)\big|dx=\varpi_j
	\end{align*}
	where $\varpi_j\to 0$, as $j\to\infty$, by the generalized Lebesgue dominated convergence.
	In a similar fashion one can treat the other terms.
	By the relative compactness of the trajectory $u(t)$ into $W^{1,p}_0(\O)$,
	there exists a subsequence $(\xi_{j_k})$,
	that we rename into $(\xi_j)$, such that $u(\xi_j)$
	is strongly convergent to some $\hat z$ in $W^{1,p}_0(\Omega)$ as $j\to\infty$.
	Then, letting $j\to\infty$ in~\eqref{quasisol},
	the generalized Lebesgue dominated convergence yields
		\begin{equation*}
	\int_{\Omega}a(\hat z)|\nabla \hat z|^{p-2}\nabla \hat z\cdot\nabla \varphi dx
	+\frac{1}{p}\int_{\Omega}a'(\hat z)|\nabla \hat z|^{p}\varphi dx \\
	=\int_{\Omega}f(\hat z)\varphi dx,\quad\forall\varphi\in C^\infty_c(\O),
	\end{equation*}
	showing that $\hat z$ is a solution of problem~\eqref{eq:pe}\footnote{Notice that we assumed  $\|\varphi\|_{C^1}\leq 1$. It is easily seen, anyway, that this assumption may be dropped via rescaling.}. Then, on one hand, we have
	$u(t_j)\to z$ in $W^{1,p}_0(\Omega)$ as $j\to\infty$ and, on the other hand,
	$u(\xi_j)\to \hat z$ in $W^{1,p}_0(\Omega)$ as $j\to\infty$. In light of
	the second part of the statement of Lemma~\ref{elle2conv}, we have
	\begin{align*}
	\|z-\hat z\|_{W^{1,p}_0(\O)} &\leq
	\|z-u(t_j)\|_{W^{1,p}_0(\O)}+
	\|u(t_j)-u(\xi_j)\|_{W^{1,p}_0(\O)}+
	\|u(\xi_j)-\hat z\|_{W^{1,p}_0(\O)} \\
	& \leq  \sup_{\mu\in [0,1]}\|u(t_j)-u(t_j+\mu)\|_{W^{1,p}_0(\Omega)}+o(1)=o(1),
	\end{align*}
	as $j\to\infty$, yielding $\hat z=z$ and concluding the proof.
\end{proof}

\begin{remark}
Forcing the nonlinearity $f$ to be zero for negative values,
the sign condition on $a'$ usually induces global solutions starting from positive
initial data to remain positive for all times $t>0$.
In fact, let us definite $\hat f:\R\to\R$ by setting
\begin{equation}
\label{nuovaL}
\hat f(s)=
\begin{cases}
f(s) & \text{if $s\geq 0$}, \\
0    & \text{if $s<0$},
\end{cases}
\end{equation}
assume that $u_0\geq 0$ a.e.\ in $\Omega$ and, furthermore, that
\begin{equation}
	\label{againsign}
a'(s)\leq 0,\quad\text{for all $s\leq 0$}.
\end{equation}
Then the solutions to the problem
\begin{equation}
	\label{problemaMod2}
\begin{cases}
u_t-{\rm div}(a(u)|\nabla u|^{p-2}\nabla u)+\frac{1}{p}a'(u)|\nabla u|^p=\hat f(u) & \text{in $(0,\infty)\times\Omega$,} \\
u(0,x)=u_0(x)  & \text{in $\Omega$,} \\
u(t,x)=0  & \text{in $(0,\infty)\times\partial\Omega$,}
\end{cases}
\end{equation}
satisfy $u(x,t)\geq 0$, for a.e.\ $x\in\Omega$ and all $t\geq  0$.
In fact, let us consider the Lipschitz function $Q:\R\to\R$ being defined by
$$
Q(s)=
\begin{cases}
0 & \text{if $s\geq 0$}, \\
s & \text{if $s\leq 0$}.
\end{cases}
$$
Testing equation~\eqref{problemaMod2} by $Q(u)\in W^{1,p}_0(\Omega)$
(which is an admissible test by~\eqref{againsign} in view
of the result of~\cite{BB} being $a'(u)|\nabla u|^pQ(u)\geq 0$ a.e. in $\R^n$)
and recalling~\eqref{nuovaL}, we get
\begin{align*}
\int_\Omega	u_tQ(u)dx+
\int_\Omega a(u)|\nabla u|^{p-2}\nabla u\nabla Q(u)dx
+\frac{1}{p}\int_\Omega a'(u)|\nabla u|^pQ(u)dx=\int_\Omega\hat f(u)Q(u)dx.
\end{align*}
Notice that it holds
$$
\int_\Omega	u_tQ(u)dx=\frac{1}{2}\frac{d}{dt}\int_\Omega	Q^2(u)dx,\qquad
\int_\Omega\hat f(u)Q(u)dx=0.
$$
as well as
\begin{align*}
& \int_\Omega a(u)|\nabla u|^{p-2}\nabla u\cdot \nabla Q(u)dx=
\int_{\Omega\cap\{u\leq 0\}} a(u)|\nabla u|^{p}dx\geq 0,   \\
&\int_\Omega a'(u)|\nabla u|^pQ(u)dx=
\int_{\Omega\cap\{u\leq 0\}} a'(u)u|\nabla u|^pdx\geq 0.
\end{align*}
In turn we conclude that
$$
\frac{d}{dt}\int_\Omega	Q^2(u(t))dx\leq 0,
$$
which yields the assertion by the definition of $Q$ and
the assumption that the initial datum $u_0$ is positive,
being $Q(u(t))=0$, for all times $t>0$.
\end{remark}

\medskip

\section{Proof of the results}
\label{prove}

Finally we can prove the main results.
\vskip6pt
\noindent
{\bf Proof of Theorem~\ref{main1}.}
Assume that $f$ is strictly positive in $(0, \infty )$ and $\O$ is strictly convex with respect
to a direction, say $x_1$, and symmetric with respect to the hyperplane $\{x_1=0\}$.
By Proposition \ref{prosimmetri2}, since $\Lambda=0$ in this case,
it follows $u(x_1,x')\leq u(-x_1,x')$ for $x_1\leq 0$.
In the same way one can prove that $u(x_1,x')\geq u(-x_1,x')$.
Therefore
$$
u(x_1,x')= u(-x_1,x'),
$$
that is $u$  belongs to the class ${\mathcal S}_{x_1}$, since
the monotonicity follows by~\eqref{zxclu:1.9222u} in Proposition~\ref{prosimmetri2}.
Finally, if $\Omega$ is a ball, by repeating this argument along any direction,
it follows that $u$ belongs to ${\mathcal R}$.

\vskip6pt
\noindent
{\bf Proof of Theorem~\ref{main-mid}.}
Part {\rm (a)} of the assertion follows by combining Theorem~\ref{main1}
with Lemma~\ref{parab-elliptic-th}. According to the notations in the statement
of Theorem~\ref{main-mid}, if $z\neq 0$ and $z\in W^{1,p}_0\cap L^\infty(\O)$
then by the regularity results of~\cite{Di,Li,T} it follows that $z\in C^{1,\alpha}(\bar{\O})$
and hence the assumptions of Theorem~\ref{main1} are fulfilled.
Part {\rm (b)} follows by combining Theorem~\ref{main1} with a
uniqueness result (of radial solutions) due to Erbe-Tang~\cite[Main Theorem, p.355]{ET}.
\vskip6pt
\noindent
{\bf Proof of Theorem~\ref{main2}.}
Part {\rm (a)} of the assertion follows from a combination of Theorem~\ref{main1}
with Lemma~\ref{lemma-bis-comp}, while part {\rm (b)} follows by combining Theorem~\ref{main1} with a
uniqueness result (of radial solutions) due to Erbe-Tang~\cite[Main Theorem, p.355]{ET}.

\bigskip

\end{document}